\documentclass[11pt,reqno]{amsart}

\usepackage[a4paper,
           top=30mm,bottom=30mm,
           left=30mm,right=30mm]{geometry}

\usepackage{amsmath,amsthm,amsfonts,amssymb,amscd,mathtools,bm}
\usepackage{xcolor}
\usepackage[unicode=true,hidelinks]{hyperref}
\usepackage[shortlabels]{enumitem}
\usepackage[T1]{fontenc}
\usepackage{empheq}
\usepackage{placeins}

\theoremstyle{plain}
\newtheorem{theorem}{Theorem}[section]
\newtheorem{lemma}[theorem]{Lemma}
\newtheorem{corollary}[theorem]{Corollary}
\newtheorem{proposition}[theorem]{Proposition}

\theoremstyle{definition}

\theoremstyle{remark}
\newtheorem*{remark}{Remark}

\numberwithin{equation}{section}
\allowdisplaybreaks

\newcommand{\dx}{\mathrm{d}x}
\newcommand{\rd}{\mathrm{d}}

\newcommand{\R}{\mathbb{R}}

\newcommand{\Z}{\mathbb{Z}}

\DeclareMathOperator{\fracpart}{frac}

\newcommand{\bee}{\begin{equation}}
\newcommand{\eee}{\end{equation}}

\title[]{\boldmath The bilinear Hilbert transforms are not of \\ restricted weak type $(1,1,1/2)$}

\author[]{Nuno J. Alves}
\author[]{Loukas Grafakos}

\address[N. J. Alves]{
CEMSE Division, King Abdullah University of Science and Technology (KAUST),
Thuwal 23955-6900, Saudi Arabia.}
\email{nuno.januarioalves@kaust.edu.sa}

\address[L. Grafakos]{
Department of Mathematics, University of Missouri, Columbia MO 65211, USA.}
\email{grafakosl@umsystem.edu}

\begin{document}

\begin{abstract}
For every real $\alpha\neq0,1$, we show that the bilinear Hilbert transform
\[
 H_\alpha(f_1,f_2)(x)
 =\text{p.v.}\int_\R f_1(x-\alpha y)\,f_2(x-y)\,\frac{\rd y}{y}
\]
fails to map $L^1(\R)\times L^1(\R)$ to $L^{1/2,\infty}(\R)$ when restricted to characteristic functions of sets of positive finite measure. The counterexample is quantitative and can be realized by sets that are finite unions of intervals of a common arbitrarily small length.
\end{abstract}

\subjclass[2020]{42B20}
\keywords{bilinear Hilbert transform, restricted weak type, endpoint estimate}
\maketitle
\thispagestyle{empty}

\section{Introduction and main result}\label{sec:introduction}
For a real parameter $\alpha$, we consider the bilinear Hilbert transform, initially defined for real-valued Schwartz functions $f_1,f_2$ on $\R$ by
\begin{equation}\label{eq:definition}
 H_\alpha(f_1,f_2)(x)
 =\text{p.v.}\int_\R
 f_1(x-\alpha y)\,f_2(x-y)\,\frac{\rd y}{y}.
\end{equation}
This family of operators arose in Calder\'on's work on commutators and the Cauchy integral on Lipschitz curves~\cite{Cal1965}, motivated by the representation of his first commutator as an average of $H_\alpha$ over the parameter $\alpha$.

When $\alpha\neq0,1$, the fundamental work of Lacey and Thiele~\cite{LT1997,LT1999} establishes
\begin{equation}\label{eq:known-bounds}
 \big\|H_\alpha(f_1,f_2)\big\|_p
 \le C_{\alpha,p_1,p_2}\,
 \|f_1\|_{p_1}\,\|f_2\|_{p_2}
\end{equation}
whenever
\[
 1<p_1,p_2\le\infty,
 \qquad
 \frac1p=\frac1{p_1}+\frac1{p_2},
 \qquad
 \frac23<p<\infty.
\]
Here, and throughout, $\|\!\cdot\!\|_p$ denotes the usual Lebesgue quasi-norm on $\R$, $0 < p \le \infty$. For an exposition of these results and their connection with Calder\'on's commutators, see Muscalu and
Schlag~\cite[Chapters~4 and~6]{MS2013}.

At the critical exponent $p=2/3$, Di Plinio~\cite{DP2026} establishes boundedness of $H_\alpha$ from $L^{p_1}(\R)\times L^{p_2}(\R)$ into $L^{2/3,\infty}(\R)$ whenever $1<p_1,p_2<2$ and $1/p_1+1/p_2=3/2$. Whether the strong-type estimate~\eqref{eq:known-bounds} holds in the remaining range $1/2<p\le2/3$, with
$1<p_1,p_2<\infty$ and $1/p=1/p_1+1/p_2$, remains open. A restricted weak-type $(1,1,1/2)$ estimate at $p_1=p_2=1$, combined with the known bounds through multilinear interpolation, would yield the estimates in this range.
Such an endpoint estimate would require a constant $C_\alpha$
such that
\begin{equation}\label{eq:restricted-weak-type}
 \lambda \,
 \big|\{x\in\R:\,
 |H_\alpha(\chi_{A_1},\chi_{A_2})(x)|>\lambda\}\big|^2
 \le C_\alpha \, |A_1| \, |A_2|
\end{equation}
for all measurable sets $A_1,A_2$ of finite measure and all $\lambda>0$.

We prove that the restricted weak-type estimate \eqref{eq:restricted-weak-type} fails for every $\alpha\in\R\setminus\{0,1\}$. More precisely, we establish the following quantitative statement. 

\begin{theorem}\label{thm:main}
For every $\alpha\in\R\setminus\{0,1\}$, there exist positive constants $c_\alpha$ and $\lambda_\alpha$ such that, for every $0<\lambda<\lambda_\alpha$ and every
$h_0>0$, there exist $0<h<h_0$ and sets $A,B\subseteq[0,1)$ of positive measure, each a finite union of half-open intervals
of length~$h$, such that
\begin{equation}\label{eq:main}
 \lambda\,
 \big|\{x\in\R:\,|H_\alpha(\chi_{A},\chi_{B})(x)|>\lambda\}\big|^2
 \ge c_\alpha\log\log\log(1/\lambda)\,|A|\,|B|.
\end{equation}
\end{theorem}

A restricted weak-type estimate with constant $C_\alpha$ as in~\eqref{eq:restricted-weak-type} would imply
\[
 c_\alpha\log\log\log(1/\lambda)\le C_\alpha
\]
for all $0<\lambda<\lambda_\alpha$. Since the
left-hand side tends to infinity as $\lambda \to 0^+$,
no such constant exists.

We now observe that the degenerate parameters $\alpha \in \{0,1\}$ have elementary descriptions. Write
\[
 Hf(x)=\text{p.v.}\int_\R
 f(x-y)\,\frac{\rd y}{y}
\]
for the linear Hilbert transform.

At $\alpha=0$, we have $H_0(f_1,f_2)=f_1Hf_2$. Given two measurable sets $A_1,A_2$ of finite measure, set
\[
 E_\lambda:=\big\{x\in\R:\,
 |H_0(\chi_{A_1},\chi_{A_2})(x)|>\lambda\big\}, \qquad \lambda> 0.
\]
Then
\[
 E_\lambda
 =A_1\cap\big\{x\in\R:\,|H\chi_{A_2}(x)|>\lambda\big\},
\]
so $|E_\lambda|\le|A_1|$, while the linear weak-type $(1,1)$ estimate gives $\lambda|E_\lambda|\le C|A_2|$. It follows that
\[
 \lambda|E_\lambda|^2\le C|A_1||A_2|.
\]

At $\alpha=1$, we have $H_1(f_1,f_2)=H(f_1f_2)$.
For $A_1=A_2=[0,1)$ and $x>1$,
\[
 H_1(\chi_{[0,1)},\chi_{[0,1)})(x)
 =\log\frac{x}{x-1}>\frac1x.
\]
Thus, for $0<\lambda<1$,
\[
 (1,\lambda^{-1})
 \subseteq
 \big\{x\in\R:\,
 |H_1(\chi_{A_1},\chi_{A_2})(x)|>\lambda\big\}.
\]
Since $|A_1||A_2|=1$, it follows that
\[
 \frac{\lambda
 \big|\{x\in\R:\,
 |H_1(\chi_{A_1},\chi_{A_2})(x)|>\lambda\}\big|^2}
 {|A_1||A_2|}
 \ge\lambda(\lambda^{-1}-1)^2
 \to\infty
 \qquad\text{as }\lambda\to0^+,
\]
so $H_1$ does not satisfy~\eqref{eq:restricted-weak-type}.

Combining these observations with Theorem~\ref{thm:main} gives immediately:

\begin{corollary}
The bilinear Hilbert transform $H_\alpha$ is of restricted weak type $(1,1,1/2)$ if and only if $\alpha=0$.
\end{corollary}

Furthermore, for $\alpha\neq0$, the change of variables $z=\alpha y$ in
\eqref{eq:definition} gives
\begin{equation}\label{eq:inversion}
 H_\alpha(f_1,f_2)(x)
 =\operatorname{sgn}(\alpha) \,
 H_{1/\alpha}(f_2,f_1)(x),
\end{equation}
and therefore, it suffices to prove
Theorem~\ref{thm:main} for
\begin{equation}\label{eq:normalized}
 \alpha \in [-1,1)\setminus \{0\}.
\end{equation}
We place this assumption on $\alpha$ throughout the rest of the manuscript, unless explicitly stated otherwise.

\smallskip

The proof of Theorem~\ref{thm:main} begins with a recursive construction of an auxiliary step function. At each step, we place enlarged copies of the preceding function on a longer interval and increase the resulting function on selected
short intervals. After dividing by a suitable constant, we approximate this function by a characteristic function. A counting argument then selects the other set while controlling cancellation in the integral. Estimates on the length of the construction interval give the triple-logarithmic lower bound.

The recursive argument used in the proof is inspired by the work of Assani and Buczolich~\cite{AB2010} on the bilinear Hardy--Littlewood maximal function and Furstenberg averages. For the bilinear Hilbert transform, Lacey and Thiele~\cite{LT1997,LT1999} established the known bounds in~\eqref{eq:known-bounds} using time-frequency analysis. These methods were adapted to the bilinear maximal function by Lacey~\cite{L2000} and to maximal multilinear averages by Demeter, Tao, and Thiele~\cite{DTT2008}. Further developments include estimates uniform in the parameter~\cite{T2002, GL2004, Li2006, UW2022}, distributional and endpoint bounds~\cite{BG,DPT}, and multilinear extrapolation methods~\cite{CGMS2009}.

We believe that the triple-logarithmic lower bound in Theorem~\ref{thm:main} is not optimal. We are working to improve this bound in a subsequent version of the manuscript. Additionally, our methods provide analogous results for the trilinear Hilbert transform and the bilinear maximal function, and may also apply to other operators. Related questions for trilinear operators were studied by Christ~\cite{C2001} and Demeter~\cite{D2008}.
Recent boundedness results for the trilinear Hilbert transform were obtained by OpenAI~\cite{OAI2026} and extended by Hu~\cite{H2026}.

The manuscript is organized as follows.
Section~\ref{sec:arithmetic} establishes an estimate for the fractional parts of arithmetic progressions. In Section~\ref{sec:integers}, we construct finite sets of integers and prove the counting estimate used in the recursive argument of Section~\ref{sec:recursive}. We then iterate the extension procedure in
Section~\ref{sec:step} to construct the auxiliary functions, which are used in Section~\ref{stripe:section} to obtain the characteristic functions. Finally, Section~\ref{sec:proof} combines the results of the preceding sections to prove Theorem~\ref{thm:main}.

\section{Fractional parts of arithmetic progressions} \label{sec:arithmetic}

In this section we prove an estimate used in
Lemma~\ref{lem:integers_source_count}.
We write $\lfloor x\rfloor$ and $\lceil x\rceil$ for the floor and ceiling of a real number $x$, respectively, and denote its fractional part by
$\fracpart(x)=x-\lfloor x\rfloor$.

\begin{lemma} \label{lem:arithmetic1}
 Let $\alpha\in(-1,1)\setminus\{0\}$ and set
\begin{equation}\label{eq:arithmetic_rho}
\rho:=\min(|\alpha|,1-|\alpha|).
\end{equation}
For every $\theta\in\mathbb R$ and every integer $m\ge1$, there exists $\theta'\in\mathbb R$ such that the two finite sequences
\[
\bigl(\fracpart(\theta+n\alpha)\bigr)_{n=0}^{m-1}
\qquad\text{and}\qquad
\bigl(\fracpart(\theta'+n\rho)\bigr)_{n=0}^{m-1}
\]
have the same terms, either in the same order or in reverse order.
\end{lemma}

\begin{proof}
 If $|\alpha|\le1/2$, then $\alpha=\rho$ or $\alpha=-\rho$. If $|\alpha|>1/2$, then $\alpha=1-\rho$ or $\alpha=-1+\rho$. Thus either $\alpha-\rho$ or $\alpha+\rho$ belongs to $\{-1,0,1\}$.

If $\alpha-\rho\in\{-1,0,1\}$, we take $\theta'=\theta$. In this case, for every integer $n$,
\[
\fracpart(\theta+n\alpha)
=\fracpart(\theta'+n\rho).
\]

Otherwise, $\alpha+\rho\in\{-1,0,1\}$ and we take
\[
\theta'=\theta-(m-1)\rho.
\]
For $0\le n<m$, we then have
\[
\fracpart(\theta+n\alpha)
=\fracpart(\theta-n\rho)
=\fracpart\bigl(\theta'+(m-1-n)\rho\bigr).
\]
As $n$ runs from $0$ to $m-1$, the index $m-1-n$ runs through the same integers in reverse order. This proves the assertion.
\end{proof}

We now bound the number of terms whose fractional parts belong to $[0,1/2)$.

\begin{lemma}\label{lem:arithmetic_estimate}
Let $\alpha\in(-1,1)\setminus\{0\}$. There exists an integer
$m_0=m_0(\alpha)\ge1$ such that, for every $\theta\in\R$ and
every integer $m\ge m_0$,
\begin{equation}\label{eq:arithmetic_estimate}
 \frac m4
 \le
 \#\Big\{n\in\Z:\, 0\le n<m,\,
 \fracpart(\theta+n\alpha)\in[0,1/2)\Big\}
 \le\frac{3m}{4}.
\end{equation}
\end{lemma}

\begin{proof}
Let $\rho$ be defined by~\eqref{eq:arithmetic_rho} and set $k=\lfloor1/(2\rho) \rfloor$.
Then $0<\rho\le1/2$ and $k\ge1$. Since $1/(2\rho)<k+1$, we have
\[
 k\rho>\frac{k}{2(k+1)}\ge\frac14.
\]
Set
\begin{equation}\label{eq:arithmetic_m_0_choice}
 m_0=m_0(\alpha):=
 \left\lceil\frac{2k}{k\rho-\frac14}\right\rceil.
\end{equation}
Then $m_0\ge1$ and
\begin{equation} \label{eq:arithmetic_m_0}
 m_0\bigg(k\rho-\frac14\bigg)\ge2k.
\end{equation}

By Lemma~\ref{lem:arithmetic1}, the number of terms whose fractional parts belong to $[0,1/2)$ is unchanged if we replace $\theta+n\alpha$ by $\theta'+n\rho$, for a suitable real number $\theta'$. We will prove the required bounds for every real
number $\theta'$. For simplicity, we write $\theta$ instead of $\theta'$ and consider
\begin{equation} \label{eq:arithmetic_finite_sequence}
 \theta,\, \theta+\rho,\ldots,\theta+(m-1)\rho.
\end{equation}

These $m$ numbers are exactly the points of $\theta+\rho\Z$ in
$[\theta,\theta+m\rho)$. Let $J$ be the number of intervals $[j,j+1)$,
$j\in\Z$, contained in $[\theta,\theta+m\rho)$. Their indices satisfy
$\lceil\theta\rceil\le j\le\lfloor\theta+m\rho\rfloor-1$, so
\[
 J=\max\bigl(0,\lfloor\theta+m\rho\rfloor-\lceil\theta\rceil\bigr)
 \ge m\rho-2.
\]
Every half-open interval of length $1/2$ contains at least $k$
points of $\theta+\rho\Z$. Indeed, after translation by $-\theta$
and division by $\rho$, this is equivalent to the fact that a half-open interval of length $1/(2\rho)$ contains at least $\lfloor1/(2\rho)\rfloor$ integers.

Therefore, in each of these $J$ intervals, both $[j,j+1/2)$ and $[j+1/2,j+1)$ contain at least $k$ points of $\theta+\rho\Z$. Since the whole interval lies in $[\theta,\theta+m\rho)$, these are points of the finite sequence in~\eqref{eq:arithmetic_finite_sequence}. Their fractional parts belong to $[0,1/2)$ and $[1/2,1)$, respectively. Since these $J$ unit intervals are disjoint, this sequence has at least $kJ$ terms with fractional part in $[0,1/2)$ and at least $kJ$ terms with fractional part in $[1/2,1)$. For $m\ge m_0$, the condition~\eqref{eq:arithmetic_m_0} gives
\[
m\left(k\rho-\frac14\right)
\ge m_0\left(k\rho-\frac14\right)
\ge 2k,
\]
and hence
\[
kJ\ge k(m\rho-2)
=\frac m4+m\left(k\rho-\frac14\right)-2k
\ge\frac m4.
\]
Thus at least $m/4$ terms have fractional part in each half of $[0,1)$.
Since the two half-open intervals partition $[0,1)$, the lower bound for the second half also gives the upper
bound $3m/4$ for the first half. This proves \eqref{eq:arithmetic_estimate}.
\end{proof}

When $\alpha=-1$, we have, for all $n \in \Z$,
\[
 \fracpart(\theta-n)=\fracpart(\theta),
\]
so \eqref{eq:arithmetic_estimate} does not hold.
We treat this case separately in the proof of
Lemma~\ref{lem:integers_source_count}.

\section{Finite sets of integers and counting} \label{sec:integers}

We construct a sequence of nested finite sets of integers and describe how their points are distributed. The main result of this section is the counting estimate in Lemma~\ref{lem:integers_source_count}, whose proof uses Lemma~\ref{lem:arithmetic_estimate} when $\alpha\ne-1$ and a separate argument when $\alpha=-1$.

\subsection{Parameters and sets}
Fix $\alpha\in[-1,1)\setminus\{0\}$. For $\alpha\ne-1$, let $m_0(\alpha)$ be as in Lemma~\ref{lem:arithmetic_estimate}, and for $\alpha=-1$, set $m_0(-1):=1$. Choose an integer $L_\alpha\ge10$ such that
\begin{equation} \label{eq:integers_L_alpha}
2^{L_\alpha-2}\ge m_0(\alpha),
\qquad
|\alpha|2^{L_\alpha-1}\ge8.
\end{equation}

Fix an integer $n\ge2$. Set $K_0=0$ and let $a_1\ge1$ be an integer. We choose the remaining parameters inductively. For $j=1,\ldots,n$, once $a_j$ has been chosen, proceed as follows:
\begin{equation}\label{eq:integers_initial_parameters}
 \begin{minipage}{0.85\linewidth}
  \begin{enumerate}[(i)]
   \item Choose an integer $\omega_j\ge a_j$.
   \item Set $K_j=\omega_j+10$.
   \item If $j<n$, choose an integer
         $a_{j+1}\ge K_j+L_\alpha$.
  \end{enumerate}
 \end{minipage}
\end{equation}
Note that this construction gives $K_1\ge11$ and implies the separation condition 
\begin{equation}\label{eq:integers_separation}
 K_j-K_{j-1} \ge L_\alpha+10 \qquad  \text{for } 2\le j \le n
\end{equation}
since, for $2\le j\le n$,
$K_j=\omega_j+10\ge a_j+10\ge K_{j-1}+L_\alpha+10$.

Now, set
\begin{equation}\label{eq:integers_Q}
 Q=2^{K_n}
\end{equation}
and call $[0,Q)$ the \textit{level-$n$ interval}. 

Starting from the level-$n$ interval $[0,Q)$, we define the remaining levels recursively. For $j=n,\ldots,1$, suppose that the level-$j$ intervals
have been defined. Divide each of them into consecutive left-closed, right-open subintervals of length $2^{K_{j-1}}$. Within each level-$j$ interval, index these subintervals from left to right by
\[
 0,1,\ldots,2^{K_j-K_{j-1}}-1.
\]
Retain the odd-indexed subintervals and call them the \textit{level-$(j-1)$ intervals}.

For each $0\le j\le n$, there are
\[
 2^{K_n-K_j-(n-j)}
\]
level-$j$ intervals, each of length $2^{K_j}$. Let $A_j$ consist of all the integer points in these intervals. Since every level-$(j-1)$ interval lies inside a level-$j$ interval, we have
\[
 A_0\subseteq A_1\subseteq\cdots\subseteq A_n
 =\{0,\ldots,Q-1\}.
\]

Set $R_1=A_1$ and, for $2\le j\le n$, define
\begin{equation} \label{eq:integers_Rj_Lj_Wj}
 R_j=A_j\setminus A_{j-1},\qquad
 L_j=2^{K_{j-1}},\qquad
 W_j=2^{K_j-8}.
\end{equation}
By construction, $A_j$ contains all the integer points in the level-$j$ intervals. For $2\le j\le n$, we obtain $A_{j-1}$ by retaining only the integer points in the odd-indexed subintervals of length $L_j$. Therefore, $R_j=A_j\setminus A_{j-1}$ consists
precisely of the integer points in the even-indexed subintervals
of the level-$j$ intervals.
The set $R_1=A_1$ contains all the integer points remaining at level $1$, including those in $A_0$. Since the sets $A_j$ are nested, the sets $R_1,\ldots,R_n$ form a partition of $\{0,\ldots,Q-1\}$. 

Finally, for each $2\le j\le n$ and each level-$j$ interval $I=[c,c+2^{K_j})$, we select the portion between $6/8$ and $7/8$ of each even-indexed subinterval of length $L_j$, and retain only the integer points that also lie in $[c+W_j,c+2^{K_j}-W_j)$.
Thus the selected set in $I$ is
\begin{equation}\label{eq:integers_good_sub}
 \mathbb Z\cap\Big[c+W_j,c+2^{K_j}-W_j\Big)\cap
 \bigcup_{\substack{0\le q<2^{K_j-K_{j-1}} \\
                   q\text{ even}}}
 \Big[c+\big(q+\tfrac68\big)L_j,\,
      c+\big(q+\tfrac78\big)L_j\Big).
\end{equation}
Let $G_j$ be the union of these selected sets over all
level-$j$ intervals.

We first compute the cardinalities of these sets.
\begin{lemma}\label{lem:integers_set_card}
For $2\le j\le n$, we have $G_j\subseteq R_j$ and
\begin{equation}\label{eq:integers_densities}
 \begin{split}
  \frac{\#A_0}{Q}&=2^{-n},\\
  r_1:=\frac{\#R_1}{Q}&=2^{-(n-1)},\qquad
  r_j:=\frac{\#R_j}{Q}=2^{-(n-j+1)},
 \end{split}
\end{equation}
as well as
\begin{equation}\label{eq:integers_good_density}
 g_j:=\frac{\#G_j}{Q}
 =\gamma2^{-(n-j)},
 \qquad
 \gamma=\frac{127}{2048}.
\end{equation}
\end{lemma}

\begin{proof}
For $j\ge2$, the integers in the even-indexed subintervals
of a level-$j$ interval belong to $R_j$. The definition of
$G_j$ therefore gives $G_j\subseteq R_j$.

At each subdivision, the number of subintervals is even,
and exactly half are retained. Since they have equal lengths
and integer endpoints, each subdivision retains half the
integers. Consequently, for $0 \le j \le n$,
\[
 \frac{\#A_j}{Q}=2^{-(n-j)}.
\]
The formula for $r_1$ follows from $R_1=A_1$. Also, for $j\ge2$, we have 
\[
 r_j
 =2^{-(n-j)}-2^{-(n-j+1)}
 =2^{-(n-j+1)}.
\]
This proves \eqref{eq:integers_densities}.

Now, fix $2\le j\le n$ and a level-$j$ interval $I=[c,c+2^{K_j})$.
By \eqref{eq:integers_separation}, both $2^{K_j}$ and~$W_j$ are multiples of
$2L_j$. Indeed,
\[
 \frac{2^{K_j}}{2L_j}=2^{K_j-K_{j-1}-1},
 \qquad
 \frac{W_j}{2L_j}=2^{K_j-K_{j-1}-9}
\]
are integers.
Thus the interval
\[
 [c+W_j,c+2^{K_j}-W_j)
\]
consists of exactly
\[
 \frac{2^{K_j}-2W_j}{2L_j}
\]
consecutive intervals of the form
$[c+2rL_j,c+2(r+1)L_j)$, with $r$ an integer.

Each of these intervals of length $2L_j$ consists of an
even-indexed subinterval of length $L_j$ followed by an
odd-indexed one. It therefore contains exactly one of the selected short intervals in~\eqref{eq:integers_good_sub}, whose length is $L_j/8$. Since $K_{j-1}\ge K_1\ge11$, this short interval has integer endpoints and contains exactly $L_j/8$ integer points. Consequently,
\[
 \#(G_j\cap I)
 =\frac{2^{K_j}-2W_j}{2L_j}\,\frac{L_j}{8}
 =\frac{2^{K_j}-2W_j}{16}
 =\left(\frac1{16}-\frac1{2^{11}}\right)2^{K_j}
 =\gamma2^{K_j}.
\]
There are $2^{K_n-K_j-(n-j)}$ level-$j$ intervals, and they
are disjoint. Summing over them gives
\[
 \#G_j
 =2^{K_n-K_j-(n-j)}\gamma\,2^{K_j}
 =\gamma \,Q\,2^{-(n-j)}.
\]
Dividing by $Q$ proves \eqref{eq:integers_good_density}.
\end{proof}

\subsection{Local structure}

We first describe how the arrangement of the integers in~$A_0$ repeats inside each level-$j$ interval.

\begin{lemma}\label{lem:integers_local_count}
Fix $2\le j\le n$ and a level-$j$ interval
$I=[c,c+2^{K_j})$.
There is a set
\[
 S\subseteq\{L_j,\ldots,2L_j-1\},
 \qquad
 \#S=L_j\,2^{-(j-1)},
\]
such that
\begin{equation} \label{eq:integers-A_0capI}
 A_0\cap I
 =
 \bigcup_{r=0}^{2^{K_j-K_{j-1}-1}-1}
 \bigl(c+2rL_j+S\bigr).
\end{equation}
The sets in this union are disjoint.
\end{lemma}

\begin{proof}
Consider an odd-indexed subinterval of $I$ of length $L_j$.
This subinterval is a level-$(j-1)$ interval with integer endpoints,
so it contains $L_j$ integer points.
To reach level~$0$, we apply the subdivision rule $j-1$ more times. At each step, we retain half the
integer points remaining from the previous step. Thus every odd-indexed subinterval contains exactly
\[
 L_j\,2^{-(j-1)}
\]
points of $A_0$.

Using the first odd-indexed subinterval of $I$, set
\[
 S:=\big\{t-c:\,t\in A_0\cap[c+L_j,c+2L_j)\big\}.
\]
Then
\[
 S\subseteq\{L_j,\ldots,2L_j-1\},
 \qquad
 \#S=L_j\,2^{-(j-1)}.
\]

Now consider a general odd-indexed subinterval of $I$,
\[
 [c+(2r+1)L_j,c+(2r+2)L_j),
\]
where $r$ is an integer with
$0\le r<2^{K_j-K_{j-1}-1}$.
This is the translate of the first odd-indexed subinterval by $2rL_j$. Since the same subdivisions are performed inside both intervals, their points belonging to $A_0$ are related by the same translation. Thus
\begin{equation}\label{eq:integers-A_0capI_aux}
 A_0\cap[c+(2r+1)L_j,c+(2r+2)L_j)
 =c+2rL_j+S.
\end{equation}
Every point of $A_0\cap I$ lies in one of these odd-indexed subintervals. Taking the union over $r$ therefore proves~\eqref{eq:integers-A_0capI}. The union is disjoint because these subintervals are disjoint.
\end{proof}

For $I$ and $S$ as in Lemma~\ref{lem:integers_local_count},
we can also describe $A_0\cap I$ by fixing $s\in S$
and varying $r$.
Since $0\le s<2L_j$, indices $r<0$ give points
below $c$, while indices
$r\ge2^{K_j-K_{j-1}-1}$ give points at or above
$c+2^{K_j}$.
Thus \eqref{eq:integers-A_0capI} can be written as
\begin{equation}\label{eq:integers-A_0capI_progressions}
 A_0\cap I
 =
 \bigcup_{s\in S}
 \Bigl((c+s+2L_j\mathbb Z)\cap I\Bigr).
\end{equation}
This union is disjoint. Indeed, an equality
\[
 c+s+2rL_j=c+s'+2r'L_j
\]
would give
\[
 s-s'=2L_j(r'-r).
\]
Since $|s-s'|<2L_j$, we must have $s=s'$ and $r=r'$.

The next lemma will be used in the case $\alpha=-1$.
For $2\le j\le n$ and a level-$j$ interval
$I=[c,c+2^{K_j})$, let
\begin{equation}\label{eq:integers_quarter_def}
 \mathcal Q_j(I):=
 \bigcup_{r=0}^{2^{K_j-K_{j-1}-1}-1}
 \Bigl[c+(2r+1)L_j,\,
       c+\bigl(2r+\tfrac54\bigr)L_j\Bigr)
\end{equation}
denote the union of the first quarters of the
odd-indexed subintervals of $I$.

\begin{lemma}\label{lem:integers_quarter_count}
Fix $2\le j\le n$ and a level-$j$ interval $I$.
For every odd-indexed subinterval $[d,d+L_j)$ of $I$,
we have
\begin{equation}\label{eq:integers_quarter_count}
 \#\bigl(A_0\cap[d,d+L_j/4)\bigr)
 =L_j\,2^{-j-1},
\end{equation}
and, for any integers $m\ge1$ and $u$ such that
$[u,u+2mL_j)\subseteq I$, 
\begin{equation}\label{eq:integers_quarter_interval_count}
 \#\bigl(A_0\cap[u,u+2mL_j)\cap\mathcal Q_j(I)\bigr)
 =m\,L_j\,2^{-j-1}.
\end{equation}
\end{lemma}

\begin{proof}
Fix an odd-indexed subinterval $[d,d+L_j)$ of $I$. At the next subdivision, its subintervals have length $2^{K_{j-2}}$. Pairing consecutive subintervals, starting at $d$, gives intervals of length
$2^{K_{j-2}+1}$. In each such pair, the even-indexed subinterval contains no
points of $A_0$. Starting from the odd-indexed subinterval, we apply the subdivision rule $j-2$ more times, retaining half the remaining integer points
at each step. Thus each pair contains exactly
\[
 2^{K_{j-2}}2^{-(j-2)}
\]
points of $A_0$.

The first quarter of $[d,d+L_j)$ is $[d,d+L_j/4)$, and
\[
 \frac{L_j/4}{2^{K_{j-2}+1}}
 =2^{K_{j-1}-K_{j-2}-3}\in\mathbb N.
\]
Indeed, for $j=2$, the exponent is $K_1-3\ge8$, since
$K_0=0$ and $K_1\ge11$. For $j\ge3$, it is larger than or equal to $17$
by \eqref{eq:integers_separation} and the fact that $L_\alpha \geq 10$.

Since the pairs start at $d$, the first quarter consists
of exactly this number of complete pairs. Therefore
\[
 \#\bigl(A_0\cap[d,d+L_j/4)\bigr)
 =
 \frac{L_j/4}{2^{K_{j-2}+1}}\,
 2^{K_{j-2}}2^{-(j-2)}
 =L_j\,2^{-j-1},
\]
which proves~\eqref{eq:integers_quarter_count}.

To prove the second assertion, write
$I=[c,c+2^{K_j})$ and let $S$ be the set given by
Lemma~\ref{lem:integers_local_count}.
Consider the first quarter
\[
 \Bigl[c+(2r+1)L_j,\,
       c+\bigl(2r+\tfrac54\bigr)L_j\Bigr)
\]
of an odd-indexed subinterval
\[
 [c+(2r+1)L_j,c+(2r+2)L_j),
\]
where $r$ is an integer with
$0\le r<2^{K_j-K_{j-1}-1}$.
By~\eqref{eq:integers-A_0capI_aux}, we have
\begin{equation} \label{eq:integers_quarter_aux}
\begin{split}
 &A_0\cap
 \Bigl[c+(2r+1)L_j,\,
       c+\bigl(2r+\tfrac54\bigr)L_j\Bigr)\\
 &\qquad=
 \bigl(c+2rL_j+S\bigr)
 \cap
 \bigl(c+2rL_j+[L_j,5L_j/4)\bigr)\\
 &\qquad=
 c+2rL_j+\bigl(S\cap[L_j,5L_j/4)\bigr).
\end{split}
\end{equation}
Putting $r=0$ and using that
translation preserves cardinality, we obtain
\[
 \#\bigl(S\cap[L_j,5L_j/4)\bigr)
 =
 \#\bigl(A_0\cap[c+L_j,c+5L_j/4)\bigr)
 =
 L_j\,2^{-j-1},
\]
where the last equality follows from
\eqref{eq:integers_quarter_count}.

Taking the union over $r$ in~\eqref{eq:integers_quarter_aux}
and using \eqref{eq:integers_quarter_def}, we obtain
\[
 A_0\cap\mathcal Q_j(I)
 =
 \bigcup_{r=0}^{2^{K_j-K_{j-1}-1}-1}
 \Bigl(c+2rL_j+\bigl(S\cap[L_j,5L_j/4)\bigr)\Bigr).
\]
Now, let $m\ge1$ and $u$ be integers such that
$[u,u+2mL_j)\subseteq I$.
Rewriting this union as in
\eqref{eq:integers-A_0capI_progressions}
and intersecting with $[u,u+2mL_j)$, we obtain
the disjoint union
\[
 A_0\cap[u,u+2mL_j)\cap\mathcal Q_j(I)
 =
 \bigcup_{s\in S\cap[L_j,5L_j/4)}
 \Bigl((c+s+2L_j\mathbb Z)\cap[u,u+2mL_j)\Bigr).
\]
For each fixed $s\in S\cap[L_j,5L_j/4)$, the condition
$c+s+2rL_j\in[u,u+2mL_j)$ can be written as
\[
 \frac{u-c-s}{2L_j}
 \le r<
 \frac{u-c-s}{2L_j}+m.
\]
This half-open interval has length $m$, so it contains exactly $m$ integers. Each set in the disjoint union therefore contains exactly $m$ points. Summing over $s\in S\cap[L_j,5L_j/4)$ and using the cardinality computed above, we obtain
\[
 \#\bigl(A_0\cap[u,u+2mL_j)\cap\mathcal Q_j(I)\bigr)
 =
 m\,\#\bigl(S\cap[L_j,5L_j/4)\bigr)
 =
 m\,L_j\,2^{-j-1},
\]
which proves~\eqref{eq:integers_quarter_interval_count}.
\end{proof}

\subsection{The counting estimate}

The proof of the next result uses Lemma~\ref{lem:arithmetic_estimate}
when $\alpha\ne-1$ and
Lemma~\ref{lem:integers_quarter_count} when $\alpha=-1$.

\begin{lemma}\label{lem:integers_source_count}
Let $2\le j\le n$, let $t_0\in G_j$, and let $k$ be
an integer with $a_j\le k\le\omega_j$.
For every $b>0$ and $0<u<b$, we have
\begin{equation}\label{eq:integers_source_conditions}
 \#\biggl\{t\in A_0:\,
 2^{k-1}\le t_0-t<2^k,\ 
 \left\lfloor t_0-\alpha(t_0-t)+\frac{u}{b}\right\rfloor
 \in R_j
 \biggr\}
 \ge 2^{k-j-3}.
\end{equation}
For each $t$ in this set, both $t$ and
$t_0-\alpha(t_0-t)+u/b$ lie in the same level-$j$ interval as $t_0$.
\end{lemma}

\begin{proof}
Fix $b>0$ and $u$ with $0<u<b$. Let $I=[c,c+2^{K_j})$ be the level-$j$ interval
containing $t_0$.
Since $t_0\in G_j$,
\begin{equation} \label{eq:integers_source_t_0_bounds}
 c+W_j\le t_0\le c+2^{K_j}-W_j-1.
\end{equation}
Recall that $W_j=2^{K_j-8}=4\cdot2^{\omega_j}$.

An integer $t$ satisfies the condition
$2^{k-1}\le t_0-t<2^k$ exactly when it belongs to
\[
 J:=[t_0-2^k+1,t_0-2^{k-1}+1).
\]
Note that this interval $J$ is contained in $I$. Indeed,
since $1\le a_j \le k\le\omega_j$,
\begin{align*}
 t_0-2^k+1
 &\ge c+W_j-2^{\omega_j}+1>c,\\
 t_0-2^{k-1}+1
 &\le t_0<c+2^{K_j}.
\end{align*}
Moreover, the interval $J$ has integer endpoints and length
$2^{k-1}=2mL_j$, where
\[
 m:=2^{k-K_{j-1}-2}
 \ge2^{L_\alpha-2}
 \ge m_0(\alpha).
\]
Here $m$ is a positive integer, and the inequalities
follow from
$k\ge a_j\ge K_{j-1}+L_\alpha$
and \eqref{eq:integers_L_alpha}.

For each integer $t\in J$, set
\[
 z(t):=t_0-\alpha(t_0-t)+\frac{u}{b}.
\]
Using $|\alpha|\le1$ and $0<u/b<1$, we obtain
\[
 |z(t)-t_0|
 <2^{\omega_j}+1<W_j.
\]
The bounds on $t_0$ in~\eqref{eq:integers_source_t_0_bounds} therefore give $z(t)\in I$.
Since the endpoints of $I$ are integers, we also have
$\lfloor z(t)\rfloor\in I$.
Thus both $t$ and $z(t)$ belong to $I$ for every integer~$t$ belonging to $J$.

Suppose first that $\alpha\ne-1$.
Let $S$ be the set given by
Lemma~\ref{lem:integers_local_count}.
Intersecting \eqref{eq:integers-A_0capI_progressions}
with $J\subseteq I$, we obtain the disjoint union
\[
 A_0\cap J
 =
 \bigcup_{s\in S}
 \Bigl((c+s+2L_j\mathbb Z)\cap J\Bigr).
\]
Since $J$ is a half-open interval of length $2mL_j$,
each set in this union contains exactly $m$ points,
spaced by $2L_j$.
For each $s\in S$, denote its smallest point by $t_s$.
Its points are therefore
\[
 t_s,\ t_s+2L_j,\ \ldots,\ t_s+2(m-1)L_j.
\]

Fix $s\in S$ and, for each $r=0,\ldots,m-1$, write
\[
 z_r:=z(t_s+2rL_j)
 =t_0-\alpha(t_0-t_s-2rL_j)+\frac{u}{b}.
\]
Then
\[
 \frac{z_r-c}{2L_j}
 =
 \frac{t_0-\alpha(t_0-t_s)+u/b-c}{2L_j}
 +r\alpha.
\]
Since $m\ge m_0(\alpha)$,
Lemma~\ref{lem:arithmetic_estimate} gives at least
$m/4$ choices of $r$ for which
\[
 \fracpart\left(\frac{z_r-c}{2L_j}\right)
 \in[0,1/2).
\]
Each $z_r$ lies in $I$, so this condition places it in the first half of one of the consecutive intervals of length $2L_j$ starting at $c$. These first halves are precisely the even-indexed subintervals of length $L_j$.
Their endpoints are integers, so $\lfloor z_r\rfloor$ lies in the same subinterval as $z_r$ and therefore belongs to $R_j$.

Thus, for each of these choices, the integer
$t_s+2rL_j$ belongs to $A_0\cap J$ and satisfies
$\lfloor z_r\rfloor\in R_j$.
It therefore belongs to the set in
\eqref{eq:integers_source_conditions}.
The candidates corresponding to different $s\in S$
are distinct, so the number of integers in that set
is at least
\[
 \frac m4\,\#S
 =
 \frac m4\,L_j2^{-(j-1)}
 =
 2^{k-j-3}.
\]

Suppose now that $\alpha=-1$.
By Lemma~\ref{lem:integers_quarter_count}, exactly
\[
 m\,L_j\,2^{-j-1}=2^{k-j-3}
\]
points of $A_0\cap J$ lie in the first quarters
of the odd-indexed subintervals of length~$L_j$. We show that every one of these points satisfies the required condition on the integer part.

Fix such a point $t$. Since $t_0,t\in\mathbb Z$ and $0<u/b<1$,
\[
 \lfloor z(t)\rfloor
 =
 \left\lfloor2t_0-t+\frac{u}{b}\right\rfloor
 =
 2t_0-t.
\]
Moreover, by the definition of $G_j$, we can write
\[
 t_0=c+(q_0+v_0)L_j,
 \qquad
 q_0\text{ even},
 \qquad
 \frac34\le v_0<\frac78,
\]
and since $t$ lies in the first quarter of an
odd-indexed subinterval, we can also write
\[
 t=c+(q+v)L_j,
 \qquad
 q\text{ odd},
 \qquad
 0\le v<\frac14.
\]
Consequently,
\[
 2t_0-t
 =
 c+(2q_0-q+1)L_j+(2v_0-v-1)L_j,
\]
where
\[
 \frac14<2v_0-v-1<\frac34.
\]
The integer $2q_0-q+1$ is even. We also have $2t_0-t=\lfloor z(t)\rfloor\in I$.
Thus the preceding expression places $2t_0-t$
in an even-indexed subinterval of $I$.
Therefore
\[
 \lfloor z(t)\rfloor=2t_0-t\in R_j.
\]
All the $2^{k-j-3}$ points counted above belong to
the set in~\eqref{eq:integers_source_conditions},
completing the proof.
\end{proof}

\section{The recursive argument}\label{sec:recursive}

In this section, we use the sets of Section~\ref{sec:integers} and Lemma~\ref{lem:integers_source_count} to prove the extension step in Proposition~\ref{prop:recur_extension}.
We will apply this step a suitable number of times in Section~\ref{sec:step} to construct the functions used in the estimates for $H_\alpha$.

Fix $\alpha\in[-1,1)\setminus\{0\}$.

\subsection{Functions and averages}

Let $p\ge1$ be an integer. At each step of the construction, we will work with a nonnegative integer-valued function $g$ on $\R$, constant on each interval $[k,k+1)$, $k\in\mathbb Z$,
and vanishing outside $[0,p)$. We also consider a function $\tau$ on $[0,p)$, constant on
the same unit intervals within $[0,p)$ and taking values in $\{2,4,8,\ldots\}$.

For each $x\in\mathbb R$, let $s_x:\mathbb R\to\mathbb R$ be given by
\begin{equation} \label{eq:recur_s_x(z)}
s_x(z)=\frac{z-(1-\alpha)x}{\alpha},
\end{equation}
and set
\begin{equation} \label{eq:recur_D_g,tau}
 D_{g,\tau}(x)
 :=\big\{z\in[0,p):\,s_x(z)<x,\quad
                  g(s_x(z)) \, \tau(z)>|z-x|\big\}.
\end{equation}
If $z\in D_{g,\tau}(x)$, then $g(s_x(z))>0$, and hence $s_x(z)\in[0,p)$. In the defining integral for~$H_\alpha$, the substitution $z=x-\alpha y$ gives $s_x(z)=x-y$. Thus the condition $s_x(z)<x$ corresponds to $y>0$.

We then define
\begin{equation}\label{eq:recur_C_g,tau}
 \begin{split}
  C_{g,\tau}(x)
  &:=\min\!\bigg(1,
       \int_{D_{g,\tau}(x)}\frac{\rd z}{\tau(z)}\bigg),\\
  \mu(g)&:=\frac{\|g\|_1}{p},\qquad
  \eta(g,\tau):=\frac1p\int_0^p C_{g,\tau}(x)\,\dx.
 \end{split}
\end{equation}
The function $(x,z) \mapsto \chi_{D_{g,\tau}(x)}(z)/\tau(z)$ is bounded and measurable on $\mathbb R\times[0,p)$, so~$C_{g,\tau}$ is measurable. The quantities $\mu(g)$ and $\eta(g,\tau)$ are the means of $g$ and $C_{g,\tau}$ on $[0,p)$. 

Our aim is to construct a pair $g,\tau$ with $\mu(g)$ arbitrarily small and $\eta(g,\tau)>1/2$.
We will start with $g=0$. In this case, the inequality
$g(s_x(z))\tau(z)>|z-x|$ cannot hold, so for every
$x\in\mathbb R$,
\[
 D_{0,\tau}(x)=\emptyset,
 \qquad C_{0,\tau}(x)=0.
\]
Thus $\mu(0)=\eta(0,\tau)=0$.
We then repeatedly apply the extension step proved below.

We first show that dilating the variable and multiplying $\tau$ by the same factor preserves $\mu(g)$ and $\eta(g,\tau)$.

\begin{lemma}\label{lem:recur_dilation}
Let $\ell\ge0$ be an integer and set $M=2^\ell$. For $p,g,\tau$ as above, the functions
\[
 g_M(x)=g(x/M),\qquad \tau_M(x)=M \,\tau(x/M)
\]
satisfy the same conditions with $Mp$ in place of $p$, and
\[
 C_{g_M,\tau_M}(Mx)=C_{g,\tau}(x),\qquad
 \mu(g_M)=\mu(g),\qquad
 \eta(g_M,\tau_M)=\eta(g,\tau).
\]
\end{lemma}

\begin{proof}
Since $M$ is a power of two, $g_M$ and $\tau_M$ satisfy the required conditions with $Mp$ in place of $p$.

Now, for every $x\in\mathbb R$ and $v\in[0,p)$, we have
\[
 s_{Mx}(Mv)=Ms_x(v),\qquad
 g_M(s_{Mx}(Mv))=g(s_x(v)),\qquad
 \tau_M(Mv)=M\tau(v).
\]
Together with $|Mv-Mx|=M|v-x|$, these identities give
\[
 D_{g_M,\tau_M}(Mx)
 =\big\{Mv:\, v\in D_{g,\tau}(x)\big\}.
\]
The change of variables $z=Mv$ therefore yields
\[
 \int_{D_{g_M,\tau_M}(Mx)}\frac{\rd z}{\tau_M(z)}
 =
 \int_{D_{g,\tau}(x)}\frac{M\,\rd v}{M\tau(v)}
 =
 \int_{D_{g,\tau}(x)}\frac{\rd v}{\tau(v)}.
\]
Taking the minimum with one proves the identity for $C_{g_M,\tau_M}$.
Finally, changing variables in the integrals defining $\mu(g_M)$ and $\eta(g_M,\tau_M)$ proves the remaining identities.
\end{proof}

\subsection{The extension step}

We will repeatedly apply the following proposition. At each step, we place dilated copies of $g$ and $\tau$ on a larger interval and increase the function obtained
from $g$ on selected unit intervals.

Let $\gamma$ be as in \eqref{eq:integers_good_density}
and, for every integer $n\ge2$, define
\begin{equation}\label{eq:recur_theta_n}
 \theta_n=\gamma |\alpha|(n-1)2^{-n-8}.
\end{equation}

\begin{proposition}\label{prop:recur_extension}
Let $p\ge1$ and $n\ge2$ be integers.
Let $g$ be a nonnegative integer-valued function on $\R$ vanishing outside $[0,p)$, and let $\tau:[0,p)\to\{2,4,8,\ldots\}$.
Assume that both functions are constant on each interval $[k,k+1)$ for $k=0,\ldots,p-1$. There exist a positive integer $P$ and functions $\widetilde g,\widetilde\tau$, satisfying the same conditions on~$[0,P)$, with
\begin{equation}\label{eq:recur_extension}
 \mu(\widetilde g)=\mu(g)+2^{1-n},\qquad
 \eta(\widetilde g,\widetilde\tau)
 \ge\eta(g,\tau)+\theta_n\big(1-\eta(g,\tau)\big),
\end{equation}
and $\|\widetilde\tau\|_\infty\le P$.

Moreover,
\begin{equation}\label{eq:recur_size_bound}
 \begin{split}
 \log_2P
 &\le2\log_2p+(4n-2)\log_2\|\tau\|_\infty\\
 &\quad+4(n-1)(L_\alpha+9)+31
 +\log_2\frac1{|\alpha|},
 \end{split}
\end{equation}
where $L_\alpha$ is as in~\eqref{eq:integers_L_alpha}.
\end{proposition}

The proof of Proposition~\ref{prop:recur_extension} will be given in two steps. First, starting from $g,\tau$, we construct the extensions $\widetilde{g},\widetilde{\tau}$ and establish all the stated bounds except the inequality for $\eta$ in~\eqref{eq:recur_extension}; this is done in Section~\ref{sec:recur_construction}. The proof is then completed in Section~\ref{sec:recur_conclusion}.

We use the following elementary approximation to determine the block size for the construction below, with numerical constants chosen for the
application in Lemma~\ref{lem:recur_common_values}.

\begin{lemma}\label{lem:recur_alignment}
Let $p$ and $Q$ be positive integers. There exist integers
$q\ge1$ and $r$ such that
\begin{equation}\label{eq:recur_alignment_b_bound}
 b:= 8pq\le\frac{2048p^2Q}{|\alpha|}
\end{equation}
and
\begin{equation}\label{eq:recur_alignment}
 |\alpha bt-8prt|\le\frac{|\alpha|}{16}
 \qquad\text{for }t=0,\ldots,Q-1.
\end{equation}
\end{lemma}

\begin{proof}
Set
\[
 R:=\bigg\lceil\frac{128pQ}{|\alpha|}\bigg\rceil
\]
and divide $[0,1)$ into $R$ half-open intervals of length $1/R$.
Two of the $R+1$ fractional parts
$\fracpart(j\alpha)$, $j=0,\ldots,R$, belong to the same interval. Let $\fracpart(j_1\alpha)$ and $\fracpart(j_2\alpha)$ be
two such fractional parts, with $0\le j_1<j_2\le R$, and set
\[
 q:=j_2-j_1,
 \qquad
 r:=\lfloor j_2\alpha\rfloor-\lfloor j_1\alpha\rfloor.
\]
Then $q,r$ are integers and $1\le q\le R$. Since the two
fractional parts lie in the same half-open interval of
length $1/R$, we have
\[
 |q\alpha-r|
 =
 \left|\fracpart(j_2\alpha)-\fracpart(j_1\alpha)\right|
 <\frac1R.
\]
Since $128pQ/|\alpha|\ge1$, we have
$R\le256pQ/|\alpha|$. Thus
\[
 b=8pq\le8pR\le\frac{2048p^2Q}{|\alpha|}.
\]
Moreover, for $0\le t<Q$,
\[
 |\alpha bt-8prt|
 =8pt\,|q\alpha-r|
 \le\frac{8pQ}{R}
 \le\frac{|\alpha|}{16}.
\]
This proves both estimates.
\end{proof}

\subsection{Construction of the extension} \label{sec:recur_construction}

Fix $p,g,\tau$ and $n$ as in Proposition~\ref{prop:recur_extension}.
In this subsection, we define $P$, $\widetilde g$, and
$\widetilde\tau$, and prove the identity for $\mu$ and the
size bounds. 

Write $\mu=\mu(g)$ and $\eta=\eta(g,\tau)$.
Since $\tau$ takes only finitely many values, all in
$\{2,4,8,\ldots\}$, there are integers $1\le a\le\omega$
such that
\[
 2^a=\min_{y\in[0,p)}\tau(y),
 \qquad
 2^\omega=\max_{y\in[0,p)}\tau(y).
\]
Set
\begin{equation} \label{eq:recur_kappa_Mj}
 \kappa:=\omega-a+L_\alpha+10,\qquad M_j:=2^{(j-1)\kappa}
\end{equation}
for every $j=1,\ldots,n$, where $L_\alpha$ is as in~\eqref{eq:integers_L_alpha}. For each such $j$, let
\begin{equation} \label{eq:recur_g^j_tau^j}
 g^{(j)}(y)=g(y/M_j),\qquad
 \tau^{(j)}(y)=M_j\,\tau(y/M_j),
\end{equation}
and define
\begin{equation} \label{eq:recur_pj}
 p_j:=M_j \, p,\qquad
 a_j:=a+(j-1)\kappa,\qquad \omega_j:=\omega+(j-1)\kappa.
\end{equation}
By Lemma~\ref{lem:recur_dilation}, these functions satisfy the
required conditions on $[0,p_j)$, and
\[
 \mu(g^{(j)})=\mu,\qquad
 \eta(g^{(j)},\tau^{(j)})=\eta.
\]
Each value of $\tau^{(j)}$ is a power of two whose exponent belongs to $[a_j,\omega_j]$.

Let $K_0=0$ and $K_j=\omega_j+10$ for $j=1,\ldots,n$.
We have $a_1=a\ge1$ and
$\omega_j-a_j=\omega-a\ge0$.
Moreover, for $1\le j<n$,
\[
 a_{j+1}=\omega_j+L_\alpha+10=K_j+L_\alpha.
\]
Thus these parameters satisfy
\eqref{eq:integers_initial_parameters}.
Form the sets $A_0,R_j,G_j$ and $Q$ as in
Section~\ref{sec:integers}.

Apply Lemma~\ref{lem:recur_alignment} with $p_n$ and this value of $Q$, and let $q,r$ be the resulting integers.
Set
\begin{equation}\label{eq:recur_b}
 b:=8p_nq.
\end{equation}
Since $p_j$ divides $p_n$ for every $j=1,\ldots,n$,
the integer $b$ is divisible by every $8p_j$ and $b\ge8$.

Let
\begin{equation}\label{eq:recur_P_J}
 P:=bQ,\qquad J:=[3b/8,5b/8),
\end{equation}
and note that $|J|=b/4$. Denote by $(g^{(j)})^\#$ and $(\tau^{(j)})^\#$ the $p_j$-periodic functions on~$\R$ that agree with $g^{(j)}$ and $\tau^{(j)}$, respectively, on $[0,p_j)$. For $t_0\in R_j$ and $0\le\xi<b$, set
\begin{equation} \label{eq:recur_g_0_tilde_tau}
 \mathfrak{g}(bt_0+\xi):=(g^{(j)})^\#(\xi),\qquad
 \widetilde\tau(bt_0+\xi):=(\tau^{(j)})^\#(\xi).
\end{equation}
These definitions specify both functions on $[0,P)$, because $R_1,\ldots,R_n$ partition $\{0,\ldots,Q-1\}$. Extend $\mathfrak{g}$ by zero outside $[0,P)$. Moreover, since $\widetilde\tau$ takes values of the periodic extensions of the functions $\tau^{(j)}$, and
$\|\tau^{(j)}\|_\infty=2^{\omega_j}$ for $j=1,\ldots,n$,
we have
\begin{equation} \label{eq:recur_tau_inf_bound}
 \|\widetilde\tau\|_\infty
 \le2^{\omega_n}
 =2^{K_n-10}
 =2^{-10}Q
 \le bQ=P.
\end{equation}

We also verify~\eqref{eq:recur_size_bound}.
By Lemma~\ref{lem:recur_alignment},
\[
 P=bQ\le\frac{2048p_n^2Q^2}{|\alpha|}.
\]
Since
\[
 p_n=p2^{(n-1)\kappa},
 \qquad
 Q=2^{K_n}=2^{\omega+(n-1)\kappa+10},
\]
we obtain
\[
 P\le\frac{p^2}{|\alpha|}
       2^{4(n-1)\kappa+2\omega+31}.
\]
Taking logarithms to base two gives
\[
 \log_2P
 \le2\log_2p+4(n-1)\kappa+2\omega
       +31+\log_2\frac1{|\alpha|}.
\]
Finally, $a\ge1$ implies
$\kappa\le\omega+L_\alpha+9$.
Substituting this bound and using
$\omega=\log_2\|\tau\|_\infty$
proves~\eqref{eq:recur_size_bound}.

We now define, for $s \in \R$,
\begin{equation}\label{eq:recur_g_tilde}
 \widetilde g(s):=\mathfrak{g}(s)+
 2b\sum_{t\in A_0}
 \chi_{[bt+b/2,\,bt+b/2+1)}(s).
\end{equation}
Since $b$ is divisible by eight and $b\ge8$, the added intervals have integer endpoints and satisfy
\[
 [bt+b/2,bt+b/2+1)
 \subseteq bt+J
 \subseteq[bt,b(t+1)).
\]
Different indices $t$ give disjoint intervals. Thus $\widetilde g$ is nonnegative and integer-valued, constant on each interval $[k,k+1)$ with $k\in\mathbb Z$, and zero outside $[0,P)$.

For $t_0\in R_j$, the definition of $\mathfrak{g}$ gives
\[
 \int_{bt_0}^{b(t_0+1)}\mathfrak{g}(s)\,\rd s
 =\int_0^b(g^{(j)})^\#(\xi)\,\rd\xi.
\]
Since $(g^{(j)})^\#$ has period $p_j$ and $b$ is a
multiple of $p_j$, we obtain
\[
 \int_0^b(g^{(j)})^\#(\xi)\,\rd\xi
 =\frac{b}{p_j}\int_0^{p_j}g^{(j)}(\xi)\,\rd\xi
 =b\mu.
\]
Each added unit interval contributes $2b$ to the integral
of $\widetilde g$. Summing over all
intervals, dividing by $P=bQ$, and using
\eqref{eq:integers_densities}, we obtain
\begin{equation}\label{eq:recur_g_tilde_mu}
 \mu(\widetilde g)
 =\sum_{j=1}^n\frac{\#R_j}{Q}\mu+
   \frac{2b\,\#A_0}{bQ}
 =\mu+2^{1-n}.
\end{equation}
It remains to prove the inequality for $\eta$
in~\eqref{eq:recur_extension}.

\subsection{Pointwise estimates}

For each $j=1,\ldots,n$, let $C_j^\#$ be the
$p_j$-periodic function on $\R$ that agrees with
$C_{g^{(j)},\tau^{(j)}}$ on $[0,p_j)$.
For $t_0\in R_j$ and $0\le\xi<b$, define
\begin{equation} \label{eq:recur_C_0}
 C_0(bt_0+\xi):=C_j^\#(\xi).
\end{equation}

\begin{lemma}\label{lem:recur_C_0_comparison}
For every $x\in[0,P)$,
\begin{equation} \label{eq:recur_C_0_comparison}
 C_{\widetilde g,\widetilde\tau}(x)
 \ge C_0(x).
\end{equation}
Moreover,
\begin{equation} \label{eq:recur_C_0_average}
 \frac1P\int_0^P C_0(x)\,\dx=\eta.
\end{equation}
\end{lemma}

\begin{proof}
For $x=bt_0+\xi$ with $t_0\in R_j$ and $0\le\xi<b$, write
$\xi=\ell p_j+x_0$, where $\ell=\lfloor\xi/p_j\rfloor$ and $0\le x_0<p_j$.
Let $c=bt_0+\ell p_j$, so that $x=c+x_0$.
The interval $c+[0,p_j)$ is contained in $[bt_0,b(t_0+1))$, because $p_j$ divides $b$.

For $z_0\in D_{g^{(j)},\tau^{(j)}}(x_0)$, both
$z_0$ and $s_{x_0}(z_0)$ belong to $[0,p_j)$.
The identity
\[
 s_{c+x_0}(c+z_0)=c+s_{x_0}(z_0)
\]
therefore shows that both $c+z_0$ and $s_{c+x_0}(c+z_0)$ belong to $c+[0,p_j)$.
On this interval, the definitions of $\widetilde\tau$
and $\widetilde g$ give
\[
 \widetilde\tau(c+z_0)=\tau^{(j)}(z_0),\qquad
 \widetilde g(c+s_{x_0}(z_0))
 \ge g^{(j)}(s_{x_0}(z_0)).
\]
It follows that
\[
 c+z_0\in
 D_{\widetilde g,\widetilde\tau}(c+x_0),
\]
and so
\[c+ D_{g^{(j)},\tau^{(j)}}(x_0) \subseteq D_{\widetilde g,\widetilde\tau}(c+x_0).\]
Using this inclusion and the change of variables $z=c+z_0$, we obtain
\begin{align*}
 \int_{D_{\widetilde g,\widetilde\tau}(x)}
       \frac{\rd z}{\widetilde\tau(z)}
 &\ge
 \int_{c+D_{g^{(j)},\tau^{(j)}}(x_0)}
       \frac{\rd z}{\widetilde\tau(z)}\\
 &=
 \int_{D_{g^{(j)},\tau^{(j)}}(x_0)}
       \frac{\rd z_0}{\widetilde\tau(c+z_0)}\\
 &=
 \int_{D_{g^{(j)},\tau^{(j)}}(x_0)}
       \frac{\rd z_0}{\tau^{(j)}(z_0)},
\end{align*}
where the last equality uses
$\widetilde\tau(c+z_0)=\tau^{(j)}(z_0)$. Since $\xi=\ell p_j+x_0$, periodicity gives
\[
 C_0(x)=C_j^\#(\xi)
       =C_j^\#(x_0)
       =C_{g^{(j)},\tau^{(j)}}(x_0).
\]
Taking the minimum with one on both sides of the integral inequality proves~\eqref{eq:recur_C_0_comparison}.

For $t_0\in R_j$, the periodicity of $C_j^\#$ and the
divisibility of $b$ by $p_j$ give
\begin{align*}
 \int_{bt_0}^{b(t_0+1)}C_0(x)\,\dx
 &=\frac b{p_j}\int_0^{p_j}
       C_{g^{(j)},\tau^{(j)}}(v)\,\rd v\\
 &=b\,\eta(g^{(j)},\tau^{(j)})
 =b \, \eta,
\end{align*}
where the last equality follows from
Lemma~\ref{lem:recur_dilation}. Since $R_1,\ldots,R_n$ form a partition of $\{0,\ldots,Q-1\}$,
summing these integrals and dividing by $P=bQ$ gives
\[
 \frac1P\int_0^P C_0(x)\,\dx
 =\sum_{j=1}^n\frac{\#R_j}{Q}\eta
 =\eta,
\]
which finishes the proof.
\end{proof}

We next introduce quantities that will be used to prove the inequality for $\eta$ in
Proposition~\ref{prop:recur_extension}.
Fix $2\le j\le n$, $t_0\in G_j$, and
$x=bt_0+\xi$ with $\xi\in J$. For each integer $k$ with $a_j\le k\le\omega_j$, set
\begin{equation} \label{eq:recur_S_k}
 \mathcal S_k
 :=\big\{t\in A_0:\,2^{k-1}\le t_0-t<2^k\big\}=A_0 \cap \big(t_0-2^k,t_0-2^{k-1}\big].
\end{equation}
The sets $\mathcal S_k$ are disjoint and
$\#\mathcal S_k\le2^{k-1}$.

For $t\in\mathcal S_k$, define
\begin{equation}\label{eq:recur_Z_t(x)}
 Z_t(x):=\int_0^1\frac{|\alpha|}{\widetilde\tau(z)}\,
       \chi_{\{\lfloor z/b\rfloor\in R_j\}}\,
       \chi_{\{\widetilde\tau(z)=2^k\}}\,\rd h,
\end{equation}
where
\[
 z=(1-\alpha)x+\alpha(bt+b/2+h).
\]
The integrand is defined to be zero when $z\notin[0,P)$,
so $\widetilde\tau$ is evaluated only on its domain.
For $t\in A_0$ belonging to none of the sets $\mathcal S_k$,
set $Z_t(x)=0$.

\begin{remark}
The assumption $t_0\in G_j$ is needed for
Lemma~\ref{lem:recur_Z_bound} below. The definitions of
$\mathcal S_k$ and $Z_t$, and Lemmas~\ref{lem:recur_Z_sets}
and~\ref{lem:recur_C-C0_bound}, remain valid under the
weaker assumption $t_0\in R_j$.
\end{remark}

\begin{lemma}\label{lem:recur_Z_sets}
Fix $2\le j\le n$, $t_0\in G_j$, and
$x=bt_0+\xi$ with $\xi\in J$.
For every $t\in A_0$,
\[
 0\le Z_t(x)\le1.
\]
Moreover, there are pairwise disjoint
measurable sets $\{E_t\}_{t \in A_0}$ such that
\[
 E_t\subseteq
 D_{\widetilde g,\widetilde\tau}(x)
 \setminus[bt_0,b(t_0+1))
\]
and
\begin{equation} \label{eq:recur_Z_t_E_t}
 Z_t(x)
 =\int_{E_t}\frac{\rd z}{\widetilde\tau(z)},
\end{equation}
for all $t \in A_0$.
\end{lemma}

\begin{proof}
For every integer $k$ with $a_j\le k\le\omega_j$ and
every $t\in\mathcal S_k$, the integrand in~\eqref{eq:recur_Z_t(x)} gives
\[
 0\le Z_t(x)\le |\alpha|2^{-k}\le1.
\]

For each integer $k$ with $a_j\le k\le\omega_j$ and each $t\in\mathcal S_k$, define the measurable set
\begin{equation} \label{eq:recur_E_t_def}
 E_t:=\Bigl\{z\in[0,P):\,
 s_x(z)\in[bt+b/2,bt+b/2+1), \ \lfloor z/b\rfloor\in R_j, \
 \widetilde\tau(z)=2^k\Bigr\}.
\end{equation}
For indices $t\in A_0$ outside all the sets $\mathcal S_k$,
set $E_t=\emptyset$.

Let $t\in\mathcal S_k$ and $z\in E_t$.
Write
\[
 s=s_x(z)=bt+b/2+h,\qquad 0\le h<1,
\]
so that $z=(1-\alpha)x+\alpha s$.
By the definition of $E_t$, we have $\widetilde\tau(z)=2^k$.
Also, $s$ belongs to the added unit interval indexed by $t$, so $\widetilde g(s)\ge2b$.
Since $\xi,b/2+h\in J$, the condition
$2^{k-1}\le t_0-t<2^k$ gives
\[
 b(2^{k-1}-1/4)<x-s<b(2^k+1/4)<2^{k+1}b.
\]
Thus $s_x(z)=s<x$, and, using $|\alpha|\le1$, we obtain
\[
 |z-x|=|\alpha|(x-s)<2^{k+1}b
 =2b\,\widetilde\tau(z)
 \le\widetilde g(s)\,\widetilde\tau(z).
\]
Together with $z\in[0,P)$, these inequalities show that
$z\in D_{\widetilde g,\widetilde\tau}(x)$.

On the other hand, using $x-s>b(2^{k-1}-1/4)$ and
$k\ge a_j\ge L_\alpha$, we obtain
\[
 |z-x|=|\alpha|(x-s)>|\alpha|b(2^{a_j-1}-1/4)>7b,
\]
where the last inequality follows from~\eqref{eq:integers_L_alpha}.
As $x\in[bt_0,b(t_0+1))$, this places $z$ outside that
interval and proves the required inclusion for $E_t$.

The sets $E_t$ are pairwise disjoint because their associated points $s_x(z)$ belong to disjoint added unit intervals. Finally, the change of variables
\[
 z=(1-\alpha)x+\alpha(bt+b/2+h)
\]
and the definition of $E_t$ prove~\eqref{eq:recur_Z_t_E_t}.
\end{proof}

\begin{lemma}\label{lem:recur_C-C0_bound}
Fix $2\le j\le n$, $t_0\in G_j$, and
$x=bt_0+\xi$ with $\xi\in J$.
We have
\begin{equation} \label{eq:recur_C-C0_bound}
 C_{\widetilde g,\widetilde\tau}(x)-C_0(x)
 \ge(1-C_0(x))
       \min\bigg(1,\sum_{t\in A_0}Z_t(x)\bigg).
\end{equation}
\end{lemma}

\begin{proof}
Write $\xi=\ell p_j+x_0$, where
$\ell=\lfloor\xi/p_j\rfloor$ and $0\le x_0<p_j$, and set $c=bt_0+\ell p_j$. Since $G_j\subseteq R_j$, the proof of Lemma~\ref{lem:recur_C_0_comparison} shows that the set
\[
 F:=c+D_{g^{(j)},\tau^{(j)}}(x_0)
\]
satisfies
\[
 F\subseteq
 D_{\widetilde g,\widetilde\tau}(x)
 \cap[bt_0,b(t_0+1))
\]
and
\[
 \int_F\frac{\rd z}{\widetilde\tau(z)}\ge C_0(x).
\]

Let $(E_t)_{t \in A_0}$ be the sets given by Lemma~\ref{lem:recur_Z_sets}.
These sets and $F$ are pairwise disjoint and contained in
$D_{\widetilde g,\widetilde\tau}(x)$.
Therefore,
\begin{align*}
 \int_{D_{\widetilde g,\widetilde\tau}(x)}
       \frac{\rd z}{\widetilde\tau(z)}
 &\ge
 \int_F\frac{\rd z}{\widetilde\tau(z)}
 +\sum_{t\in A_0}\int_{E_t}\frac{\rd z}{\widetilde\tau(z)}\\
 &\ge C_0(x)+\sum_{t\in A_0}Z_t(x).
\end{align*}
Taking the minimum with one gives
\[
 C_{\widetilde g,\widetilde\tau}(x)
 \ge\min\bigg(1,C_0(x)
                   +\sum_{t\in A_0}Z_t(x)\bigg).
\]
Finally, for $0\le v\le1$ and $w\ge0$,
\[
 \min(1,v+w)-v
 =\min(1-v,w)
 \ge(1-v)\min(1,w).
\]
Applying this inequality with $v=C_0(x)$ and
$w=\sum_{t\in A_0}Z_t(x)$ proves the required estimate.
\end{proof}

\subsection{Contribution of the added intervals}

We first use Lemma~\ref{lem:recur_alignment} to obtain a common set
on which the relevant periodic values agree. We then apply
Lemma~\ref{lem:integers_source_count} to bound the contribution
of the added intervals.

Let $b,Q,r,p_1,\ldots,p_n$ and the periodic functions
$(\tau^{(j)})^\#$, $j=1,\ldots,n$, be as in Section~\ref{sec:recur_construction}.
For $x\in\R$, $0\le h<1$, and $t=0,\ldots,Q-1$, define
\begin{equation} \label{eq:recur_zeta_t_z_t}
 \zeta_x(h):=(1-\alpha)x+\alpha(b/2+h),
 \qquad
 z_t(h):=\zeta_x(h)+\alpha bt.
\end{equation}

\begin{lemma}\label{lem:recur_common_values}
Fix $x\in\R$ and let
\[
 F_x:=\left\{h\in[0,1):\,
       \operatorname{dist}(\zeta_x(h),\mathbb Z)>
       \frac{|\alpha|}{16}\right\}.
\]
For every $j=1,\ldots,n$, every $t=0,\ldots,Q-1$,
and every $h\in F_x$, we have
\begin{equation}\label{eq:recur_common_value}
 (\tau^{(j)})^\#(z_t(h))
 =(\tau^{(j)})^\#(\zeta_x(h)).
\end{equation}
Moreover,
\begin{equation}\label{eq:recur_good_h_measure}
 |F_x|\ge\frac34.
\end{equation}
\end{lemma}

\begin{proof}
Set $\varepsilon:=|\alpha|/16$. By Lemma~\ref{lem:recur_alignment} we have
\begin{equation} \label{eq:recur_aux_align}
 |z_t(h)-8p_nrt-\zeta_x(h)|
 =|\alpha bt-8p_nrt|
 \le\varepsilon
\end{equation}
for every $t=0,\ldots,Q-1$ and $0\le h<1$.

Fix $j\in\{1,\ldots,n\}$ and $t\in\{0,\ldots,Q-1\}$.
Since $p_j$ divides $p_n$ and $r,t$ are integers,
\[
 8p_nrt=p_j\left(8rt\,\frac{p_n}{p_j}\right)
\]
is an integer multiple of $p_j$.
The $p_j$-periodicity of $(\tau^{(j)})^\#$ therefore gives
\[
 (\tau^{(j)})^\#(z_t(h))
 =
 (\tau^{(j)})^\#(z_t(h)-8p_nrt).
\]

Now let $h\in F_x$ and set
$\ell:=\lfloor\zeta_x(h)\rfloor$.
The definition of $F_x$ gives
\[
 \ell+\varepsilon<\zeta_x(h)<\ell+1-\varepsilon.
\]
Combining these inequalities with~\eqref{eq:recur_aux_align} we obtain
\[
 \ell
 <\zeta_x(h)-\varepsilon
 \le z_t(h)-8p_nrt
 \le\zeta_x(h)+\varepsilon
 <\ell+1.
\]
Thus $z_t(h)-8p_nrt$ and $\zeta_x(h)$ belong to the same
unit interval $[\ell,\ell+1)$.
Since $(\tau^{(j)})^\#$ is constant on this interval,
\[
 (\tau^{(j)})^\#(z_t(h)-8p_nrt)
 =
 (\tau^{(j)})^\#(\zeta_x(h)).
\]
Together with periodicity, this proves
\eqref{eq:recur_common_value}.

It remains to estimate $|F_x|$.
Note that $\zeta_x([0,1))$ is an interval of length $|\alpha|\le1$ that intersects at most two intervals of the form
\[
 [v-\varepsilon,v+\varepsilon],
 \qquad v\in\mathbb Z.
\]
Indeed, an interval intersecting three such neighborhoods would have length at least
\[
 2-2\varepsilon\ge\frac{15}{8}>1.
\]

Since $\zeta_x(h)=\zeta_x(0)+\alpha h$, for each integer $v$
the condition
\[
 \zeta_x(h)\in[v-\varepsilon,v+\varepsilon]
\]
is equivalent to
\[
 \left|h-\frac{v-\zeta_x(0)}{\alpha}\right|
 \le\frac{\varepsilon}{|\alpha|}.
\]
Thus the set
\[
 \bigl\{h\in[0,1):\,
       \zeta_x(h)\in[v-\varepsilon,v+\varepsilon]\bigr\}
\]
is the intersection of $[0,1)$ with an interval of length
$2\varepsilon/|\alpha|$, and therefore has measure at most
$2\varepsilon/|\alpha|$.

By the definition of $F_x$,
\[
 [0,1)\setminus F_x
 =
 \bigcup_{v\in\mathbb Z}
 \bigl\{h\in[0,1):\,
       \zeta_x(h)\in[v-\varepsilon,v+\varepsilon]\bigr\}.
\]
Since $\zeta_x([0,1))$ intersects at most two of these integer neighborhoods, at most two sets in this union are nonempty. Consequently,
\[
 |[0,1)\setminus F_x|
 \le \frac{4\varepsilon}{|\alpha|}.
\]
It follows that
\[
 |F_x|
 =1-|[0,1)\setminus F_x|
 \ge1-\frac{4\varepsilon}{|\alpha|}
 =\frac34.
\]
This proves~\eqref{eq:recur_good_h_measure}.
\end{proof}

Now, recall the quantity $Z_t(x)$ defined in~\eqref{eq:recur_Z_t(x)}.

\begin{lemma}\label{lem:recur_Z_bound}
Fix $2\le j\le n$, $t_0\in G_j$, and
$x=bt_0+\xi$ with $\xi\in J$. Then
\begin{equation}\label{eq:recur_Z_bound}
 \sum_{t\in A_0}Z_t(x)\ge|\alpha|2^{-j-4}.
\end{equation}
\end{lemma}

\begin{proof}
For this $x$, let $F_x$ be the set defined in
Lemma~\ref{lem:recur_common_values}.
Thus $|F_x|\ge3/4$, and~\eqref{eq:recur_common_value}
holds for every $h\in F_x$ and every $t\in A_0$.

For $t\in A_0$ and $0\le h<1$, set
\[
 V_t(h):=\chi_{\{\lfloor z_t(h)/b\rfloor\in R_j\}}.
\]
If $V_t(h)=1$, then $\ell:=\lfloor z_t(h)/b\rfloor \in R_j$ and
$z_t(h)\in[b\ell,b(\ell+1))\subseteq[0,P)$.
By the definition of $\widetilde\tau$ and the divisibility of $b$ by $p_j$,
\[
 \widetilde\tau(z_t(h))
 =(\tau^{(j)})^\#(z_t(h)-b\ell)
 =(\tau^{(j)})^\#(z_t(h)).
\]
When $V_t(h)=0$, the integrand defining $Z_t(x)$ vanishes. Therefore
\[
 Z_t(x)=|\alpha|2^{-k}\int_0^1
 V_t(h)\,
 \chi_{\{(\tau^{(j)})^\#(z_t(h))=2^k\}}\,\rd h,
 \qquad t\in\mathcal S_k.
\]

For each $0\le h<1$, let
\[
 e(h):=\frac{(1-\alpha)\xi+\alpha(b/2+h)}b.
\]
Since $b\ge8$, both $\xi$ and $b/2+h$ belong to
$J=[3b/8,5b/8)$, whose length is $b/4$.
Using $|\alpha|\le1$, we obtain
\[
 \left|e(h)-\frac{\xi}{b}\right|
 =|\alpha|\frac{|b/2+h-\xi|}{b}
 <\frac14.
\]
Together with $3/8\le\xi/b<5/8$, this gives
\[
 \frac18<e(h)<\frac78.
\]
In particular, $0<be(h)<b$.

Using $x=bt_0+\xi$, we also have
\[
 \frac{z_t(h)}b=t_0-\alpha(t_0-t)+e(h).
\]
For every integer $k$ with $a_j\le k\le\omega_j$,
the set $\mathcal S_k$ consists precisely of the
$t\in A_0$ satisfying
\[
 2^{k-1}\le t_0-t<2^k.
\]
We may therefore apply Lemma~\ref{lem:integers_source_count}
with its parameter $u$ equal to $be(h)$.
Its integer-part condition is exactly $V_t(h)=1$, so
\begin{equation}\label{eq:recur_source_sum}
 \sum_{t\in\mathcal S_k}V_t(h)\ge2^{k-j-3}.
\end{equation}

The sets $\mathcal S_k$ are pairwise disjoint, and $Z_t(x)=0$
for indices outside their union.
Since all sums are finite and all integrands are nonnegative,
we may group the terms by $k$ and restrict the integrals to $F_x$.
Using~\eqref{eq:recur_common_value} and then
\eqref{eq:recur_source_sum}, we obtain
\begin{align*}
 \sum_{t\in A_0}Z_t(x)
 &\ge|\alpha|\sum_{k=a_j}^{\omega_j}2^{-k}
   \int_{F_x}
   \chi_{\{(\tau^{(j)})^\#(\zeta_x(h))=2^k\}}
   \sum_{t\in\mathcal S_k}V_t(h)\,\rd h\\
 &\ge|\alpha|2^{-j-3}
   \int_{F_x}\sum_{k=a_j}^{\omega_j}
   \chi_{\{(\tau^{(j)})^\#(\zeta_x(h))=2^k\}}\,\rd h\\
 &=|\alpha|2^{-j-3}|F_x|\\
 &\ge\frac34|\alpha|2^{-j-3}
 \ge|\alpha|2^{-j-4}.
\end{align*}
The equality holds because every value of
$(\tau^{(j)})^\#$ equals $2^k$ for exactly one integer
$k\in[a_j,\omega_j]$, and the bound for $|F_x|$ follows
from~\eqref{eq:recur_good_h_measure}.
This proves~\eqref{eq:recur_Z_bound}.
\end{proof}

\subsection{Completing the extension} \label{sec:recur_conclusion}

We have defined $P$, $\widetilde g$, and $\widetilde\tau$
and verified the identity for $\mu$ and the size bounds.
We now prove the inequality for $\eta$.

\begin{proof}[Proof of Proposition~\ref{prop:recur_extension}]
Fix $2\le j\le n$, $t_0\in G_j$, and
$x=bt_0+\xi$ with $\xi\in J$.
Lemma~\ref{lem:recur_Z_bound} and the inequality
$|\alpha|2^{-j-4}\le1$ give
\[
 \min\bigg(1,\sum_{t\in A_0}Z_t(x)\bigg)
 \ge|\alpha|2^{-j-4}.
\]
Combining this bound with
Lemma~\ref{lem:recur_C-C0_bound}, we obtain
\begin{equation}\label{eq:recur_extension_pointwise}
 C_{\widetilde g,\widetilde\tau}(x)-C_0(x)
 \ge(1-C_0(x))|\alpha|2^{-j-4}.
\end{equation}
By Lemma~\ref{lem:recur_C_0_comparison}, the difference
on the left is nonnegative for every $x\in[0,P)$.

For $t_0\in G_j\subseteq R_j$, the definition of $C_0$ gives
\[
 C_0(bt_0+\xi)=C_j^\#(\xi).
\]
Since $b$ is divisible by $8p_j$, both endpoints of $J=[3b/8,5b/8)$ are multiples of $p_j$. Thus
\begin{align*}
 \int_{bt_0+J}(1-C_0(x))\,\dx
 &=\int_J(1-C_j^\#(\xi))\,\rd\xi\\
 &=\frac{|J|}{p_j}\int_0^{p_j}
       \bigl(1-C_{g^{(j)},\tau^{(j)}}(v)\bigr)\,\rd v\\
 &=\frac b4(1-\eta).
\end{align*}
The sets $bt_0+J$, with $t_0\in G_j$ and $2\le j\le n$,
are disjoint. Indeed, each is contained in
$[bt_0,b(t_0+1))$, and the indices $t_0$ are distinct
because $G_j\subseteq R_j$ and $R_1,\ldots,R_n$
form a partition of $\{0,\ldots,Q-1\}$.

Integrate~\eqref{eq:recur_extension_pointwise} over each interval $bt_0+J$, sum over $t_0\in G_j$ and $2\le j\le n$, and divide by $P=bQ$.
Using~\eqref{eq:recur_C_0_comparison}--\eqref{eq:recur_C_0_average} and
\eqref{eq:integers_good_density}, we obtain
\begin{align*}
 \eta(\widetilde g,\widetilde\tau)
 &\ge\eta+\frac{|\alpha|}{64}
   \sum_{j=2}^n\frac{\#G_j}{Q}2^{-j}(1-\eta)\\
 &=\eta+\frac{\gamma|\alpha|}{64}
   (n-1)2^{-n}(1-\eta)\\
 &=\eta+4\theta_n(1-\eta)\\
 &\ge\eta+\theta_n(1-\eta),
\end{align*}
where the last inequality uses $0\le\eta\le1$.
Together with~\eqref{eq:recur_g_tilde_mu},
\eqref{eq:recur_tau_inf_bound}, and \eqref{eq:recur_size_bound}, this completes the proof of the proposition.
\end{proof}

\section{Auxiliary step functions} \label{sec:step}

We first establish an integral representation for $H_\alpha$ and a
finite counting estimate. We then construct functions $g_0,\tau_0$ by iterating
Proposition~\ref{prop:recur_extension} and prove the integral and size
bounds needed for the subsequent argument. Finally, we establish
subdivision estimates that will be used in the next section to
pass to characteristic functions.

\subsection{Representation and counting}

For a positive integer $P$ and a function~$g$ vanishing outside $[0,P)$, we set
\begin{equation}\label{eq:step_X_set}
 X_g:=\bigg\{x\in[0,P):\,
       \operatorname{dist}(x,\operatorname{supp}g)
       \ge\frac12\bigg\}.
\end{equation}

We first provide an integral representation for $H_\alpha$
on $X_g$.

\begin{lemma}\label{lem:step_absolute_convergence}
Let $P$ be a positive integer and let $g\in L^1(\R)$ be
a nonnegative function vanishing outside $[0,P)$.
For every $f\in L^\infty(\R)$ and every $x\in X_g$,
the integral defining $H_\alpha(f,g)(x)$ converges absolutely
and
\begin{equation}\label{eq:step_kernel_representation}
 H_\alpha(f,g)(x)
 =\operatorname{sgn}(\alpha)
  \int_\R\frac{f(z)\,g(s_x(z))}{x-z}\,\rd z,
\end{equation}
where $s_x$ is defined in~\eqref{eq:recur_s_x(z)}.
\end{lemma}

\begin{proof}
Let $f\in L^\infty(\R)$ and $x\in X_g$.
If $|y|<1/2$, then
\[
 |x-(x-y)|=|y|<\frac12
 \le\operatorname{dist}(x,\operatorname{supp}g),
\]
so $x-y\notin\operatorname{supp}g$.
Hence $g(x-y)=0$ for almost every $y$ with $|y|<1/2$,
and so, for every $0<\varepsilon<1/2$,
\begin{align*}
 \int_{|y|>\varepsilon}
 \frac{|f(x-\alpha y)\,g(x-y)|}{|y|}\,\rd y
 &=
 \int_{|y|\ge1/2}
 \frac{|f(x-\alpha y)\,g(x-y)|}{|y|}\,\rd y\\
 &\le
 2\|f\|_\infty\int_\R |g(x-y)|\,\rd y\\
 &=2\|f\|_\infty\|g\|_1<\infty.
\end{align*}
Thus the defining integral converges absolutely for $x \in X_g$, and in this case we have
\begin{align*}
 H_\alpha(f,g)(x)
 &=
 \lim_{\varepsilon\to0^+}
 \int_{|y|>\varepsilon}
 \frac{f(x-\alpha y)\,g(x-y)}{y}\,\rd y\\
 &=
 \int_\R
 \frac{f(x-\alpha y)\,g(x-y)}{y}\,\rd y.
\end{align*}
Making the change of variables
$z=x-\alpha y$ gives
\[
 y=\frac{x-z}{\alpha},
 \qquad
 x-y=s_x(z),
\]
and we deduce
\begin{align*}
 H_\alpha(f,g)(x)
 &=
 \int_\R
 f(z)\,g(s_x(z))\,\frac{\alpha}{x-z}\,
 \frac{\rd z}{|\alpha|}\\
 &=
 \operatorname{sgn}(\alpha)
 \int_\R\frac{f(z)\,g(s_x(z))}{x-z}\,\rd z,
\end{align*}
which proves~\eqref{eq:step_kernel_representation}.
\end{proof}

We next prove a finite counting estimate that will be used to choose the first input.

\begin{lemma}\label{lem:step_counting}
Let $q\ge1$ and $w_1,\ldots,w_q\ge2$ be integers,
and let $a_1,\ldots,a_q$ be real numbers.
Let $t>0$ and $J_0\subseteq\{1,\ldots,q\}$.
Suppose that $a_i>2t$ for every $i\in J_0$, and set
\[
 W:=\prod_{i=1}^q\{0,\ldots,w_i-1\},\qquad
 V_s:=\sum_{i=1}^q a_i\,\chi_{\{0\}}(s_i)
\]
for $s=(s_1,\ldots,s_q)\in W$. Then
\begin{equation}\label{eq:step_counting}
 \#\big\{s\in W:\,|V_s|>t\big\}
 \ge\frac{\#W}{4}
 \min\!\bigg(1,\sum_{i\in J_0}\frac1{w_i}\bigg).
\end{equation}
\end{lemma}

\begin{proof}
If $J_0$ is empty, the right-hand side of
\eqref{eq:step_counting} is zero, so there is nothing to prove.
Otherwise, set
\[
 Y:=\big\{s\in W:\,s_i\in\{0,1\}
       \text{ for some }i\in J_0\big\}.
\]
For each $s\in Y$, let $i(s)$ be the smallest index in $J_0$ for which $s_{i(s)}\in\{0,1\}$.
Define $s'\in W$ by changing this coordinate from zero to one or from one to zero, and leaving all other coordinates unchanged. The smallest index with the stated property is the same for $s'$ as for $s$, so applying this operation to $s'$ returns $s$.
Thus it partitions $Y$ into disjoint pairs $\{s,s'\}$.

For each such pair,
\[
 |V_s-V_{s'}|=a_{i(s)}>2t.
\]
At least one of $|V_s|$ and $|V_{s'}|$ therefore exceeds
$t$, since otherwise the triangle inequality would give
$|V_s-V_{s'}|\le2t$.
Consequently,
\[
 \#\big\{s\in W:\,|V_s|>t\big\}\ge\frac{\#Y}{2}.
\]

An element $s\in W$ does not belong to $Y$ precisely when $s_i$ is neither zero nor one for every $i\in J_0$.
Counting these elements and subtracting their number from $\#W$ gives
\[
 \#Y=\#W\bigg(1-\prod_{i\in J_0}
                   \bigg(1-\frac2{w_i}\bigg)\bigg).
\]
Set $q_0:=\sum_{i\in J_0}1/w_i$.
Since $0\le2/w_i\le1$, we have
\[
 \prod_{i\in J_0}\bigg(1-\frac2{w_i}\bigg)
 \le
 \frac1{\displaystyle\prod_{i\in J_0}
                     \bigg(1+\frac2{w_i}\bigg)}
 \le\frac1{1+2q_0}.
\]
Here the first inequality follows from
$(1-u)(1+u)\le1$ for $0\le u\le1$, and the second follows by expanding the product in the denominator.
It follows that
\begin{align*}
 \#\big\{s\in W:\,|V_s|>t\big\}
 &\ge\frac{\#W}{2}\left(1-\frac1{1+2q_0}\right)\\
 &=\#W\,\frac{q_0}{1+2q_0}\\
 &\ge\frac{\#W}{4}\min(1,q_0),
\end{align*}
where in the last inequality we used
$1+2q_0\le4\max(1,q_0)$. This proves
\eqref{eq:step_counting}.
\end{proof}

\subsection{The initial functions}
\label{subsec:step_initial_functions}

Now, we construct $g_0,\tau_0$ by iterating the extension step of
Proposition~\ref{prop:recur_extension}. Recall the function
$C_{g,\tau}$ defined in~\eqref{eq:recur_C_g,tau}.

Fix an integer $N\ge8$, set
\begin{equation}\label{eq:step_A_alpha}
 A_\alpha:=\frac{512}{\gamma|\alpha|},
\end{equation}
where $\gamma$ is as in~\eqref{eq:integers_good_density}, and let
\begin{equation}\label{eq:step_n_d}
 n:=\lceil N(A_\alpha+1)\rceil+1,
 \qquad
 d:=\lceil1/\theta_n\rceil,
\end{equation}
where $\theta_n$ is defined in~\eqref{eq:recur_theta_n}.

Set $p_0:=1$, let $g^{[0]}$ be the zero function on $\R$,
and let $\tau^{[0]}$ be the constant function equal to two
on $[0,1)$. For each $i=0,\ldots,d-1$, apply
Proposition~\ref{prop:recur_extension} to
$p_i,g^{[i]},\tau^{[i]}$ and denote the resulting integer
and functions by $p_{i+1},g^{[i+1]},\tau^{[i+1]}$.

Set
\begin{equation}\label{eq:step_initial_data}
 P_0:=p_d,\qquad
 g_0:=g^{[d]},\qquad
 \tau_0:=\tau^{[d]},
\end{equation}
and let $X_0:=X_{g_0}$, where $X_{g_0}$ is defined
by~\eqref{eq:step_X_set} with $P=P_0$.

\begin{lemma}\label{lem:step_Cgtau_bound}
The number $P_0$ is a positive integer, $g_0$ is a
nonnegative integer-valued function on $\R$ vanishing
outside $[0,P_0)$, and
$\tau_0:[0,P_0)\to\{2,4,8,\ldots\}$.
Both functions are constant on each interval $[k,k+1)$
for $k=0,\ldots,P_0-1$.
Moreover,
\begin{equation}\label{eq:step_g1_bound}
 0<\|g_0\|_1\le\frac{P_0}{N},
\end{equation}
and
\begin{equation}\label{eq:step_Cgtau_bound}
 \int_{X_0}C_{g_0,\tau_0}(x)\,\dx>\frac{3P_0}{8}.
\end{equation}
\end{lemma}

\begin{proof}
The stated properties of $P_0,g_0,\tau_0$ follow from
Proposition~\ref{prop:recur_extension}.
Since
\[
 \mu(g^{[0]})=0,
 \qquad
 \eta(g^{[0]},\tau^{[0]})=0,
\]
the two estimates in~\eqref{eq:recur_extension} give,
by induction,
\begin{equation}\label{eq:step_Cgtau_aux_induction}
 \mu(g^{[i]})=i2^{1-n},
 \qquad
 \eta(g^{[i]},\tau^{[i]})
 \ge1-(1-\theta_n)^i,
 \qquad i=0,\ldots,d.
\end{equation}

Using~\eqref{eq:step_Cgtau_aux_induction} with $i=d$,
the bound $d\le1/\theta_n+1$, and the definitions of
$\theta_n$ and $A_\alpha$, we obtain
\[
 \frac{\|g_0\|_1}{P_0}
 =\mu(g^{[d]})
 =d2^{1-n}
 \le\frac{A_\alpha}{n-1}+2^{1-n}
 \le\frac{A_\alpha+1}{n-1}
 \le\frac1N.
\]
Here we used $2^{1-n}\le1/(n-1)$ for $n\ge2$ and
$n-1\ge N(A_\alpha+1)$.
This proves~\eqref{eq:step_g1_bound}.

Moreover, \eqref{eq:step_Cgtau_aux_induction} with $i=d$,
together with $0<\theta_n<1$ and $d\theta_n\ge1$, gives
\begin{align*}
 \frac1{P_0}\int_0^{P_0}C_{g_0,\tau_0}(x)\,\dx
 &=\eta(g^{[d]},\tau^{[d]})\\
 &\ge1-(1-\theta_n)^d
 >1-e^{-d\theta_n}
 \ge1-e^{-1}
 >\frac58,
\end{align*}
where in the second inequality we used the fact that
for $0<t<1$ we have
\[
 \log(1-t)
 =-\int_0^t\frac{\rd s}{1-s}
 <-\int_0^t1\,\rd s
 =-t.
\]

Now, let $N_0$ be the number of unit intervals on which
$g_0$ is positive.
Since $g_0$ is integer-valued, it is at least one on
each of these intervals. Thus
\[
 N_0\le\|g_0\|_1.
\]
For each unit interval $[k,k+1)$ on which $g_0$ is
positive, consider the interval $(k-1/2,k+3/2)$.
By the definition of $X_0$, these $N_0$ enlarged
intervals cover $[0,P_0)\setminus X_0$.
Since each has length two,
\[
 \big|[0,P_0)\setminus X_0\big|
 \le2N_0
 \le2\|g_0\|_1
 \le\frac{2P_0}{N}.
\]
Using $0\le C_{g_0,\tau_0}\le1$ and $N\ge8$, we obtain
\[
 \int_{X_0}C_{g_0,\tau_0}(x)\,\dx
 >\frac{5P_0}{8}-\frac{2P_0}{N}
 \ge\frac{3P_0}{8}.
\]
This proves~\eqref{eq:step_Cgtau_bound}.
\end{proof}

\subsection{Size bounds}

In the next two lemmas we bound $P_0$ and
$\|\tau_0\|_\infty$.
For $i=0,\ldots,d$, set
\[
 \beta_i:=\log_2p_i.
\]

\begin{lemma}\label{lem:step_size_recurrence}
For $i=1,\ldots,d$, we have
\begin{equation}\label{eq:step_target_height}
 \|\tau^{[i]}\|_\infty\le p_i.
\end{equation}
Moreover, for $i=0,\ldots,d-1$,
\begin{equation}\label{eq:step_size_recurrence}
 \beta_{i+1}
 \le4n\beta_i+4(n-1)(L_\alpha+10)+33
       +\log_2\frac1{|\alpha|}.
\end{equation}
\end{lemma}

\begin{proof}
For $i=1,\ldots,d$, applying
Proposition~\ref{prop:recur_extension} to
$p_{i-1},g^{[i-1]},\tau^{[i-1]}$
gives~\eqref{eq:step_target_height}.
Let $\omega_i$ denote the largest exponent occurring
in $\tau^{[i]}$. Then
\[
 \omega_i\le\beta_i
 \qquad\text{for }i=1,\ldots,d.
\]
Since $\omega_0=1$ and $\beta_0=0$, it follows that
$\omega_i\le\beta_i+1$ for $i=0,\ldots,d$.
Thus, for $i=0,\ldots,d-1$,
\eqref{eq:recur_size_bound} gives
\begin{align*}
 \beta_{i+1}
 &\le2\beta_i+(4n-2)\omega_i
       +4(n-1)(L_\alpha+9)+31
       +\log_2\frac1{|\alpha|}\\
 &\le2\beta_i+(4n-2)(\beta_i+1)
       +4(n-1)(L_\alpha+9)+31
       +\log_2\frac1{|\alpha|}\\
 &=4n\beta_i+4(n-1)(L_\alpha+10)+33
       +\log_2\frac1{|\alpha|},
\end{align*}
which proves~\eqref{eq:step_size_recurrence}.
\end{proof}

We now estimate $P_0=p_d$ and the $L^\infty$ norm of $\tau_0=\tau^{[d]}$.
Set
\begin{equation}\label{eq:step_T_alpha_C_alpha}
 T_\alpha:=L_\alpha+10+\log_2\frac1{|\alpha|},
 \qquad
 B_\alpha:=A_\alpha+T_\alpha+1.
\end{equation}

\begin{lemma}\label{lem:step_P_0_tau_0_size}
We have
\begin{equation}\label{eq:step_P_0_tau_0_size}
 \log(16P_0^2)\le\exp(B_\alpha2^n),
 \qquad
 \|\tau_0\|_\infty\le P_0.
\end{equation}
\end{lemma}

\begin{proof}
Since $d\ge1$, \eqref{eq:step_target_height} gives
\[
 \|\tau_0\|_\infty
 =\|\tau^{[d]}\|_\infty
 \le p_d=P_0.
\]

To estimate $P_0$, note that
\begin{align*}
 &(4n-1)T_\alpha-4(n-1)(L_\alpha+10)-33
       -\log_2\frac1{|\alpha|}\\
 &\qquad=3L_\alpha-3
       +(4n-2)\log_2\frac1{|\alpha|}
 \ge0.
\end{align*}
Here we used $L_\alpha\ge10$ and $|\alpha|\le1$.
Consequently, \eqref{eq:step_size_recurrence} implies
\[
 \beta_{i+1}+T_\alpha
 \le4n(\beta_i+T_\alpha),
 \qquad i=0,\ldots,d-1.
\]
Since $\beta_0=0$, induction gives
\begin{equation}\label{eq:step_size_iteration}
 \beta_d\le T_\alpha\bigl((4n)^d-1\bigr).
\end{equation}

The definitions of $d$ and $A_\alpha$ give
\[
 d\le\frac{A_\alpha}{2}\frac{2^n}{n-1}+1.
\]
Our choice of $n$ in~\eqref{eq:step_n_d} gives $n\ge3$.
For these integers,
\[
 \log(4n)\le2(n-1),
 \qquad
 \log(4n)\le2^n.
\]
Indeed, $\log n\le n-1$ and $\log4\le n-1$ prove
the first inequality, and $2(n-1)\le2^n$ proves the second.
Therefore,
\begin{align*}
 d\log(4n)
 &\le
 \bigg(\frac{A_\alpha}{2}\frac{2^n}{n-1}+1\bigg)
       \log(4n)\\
 &\le A_\alpha2^n+\log(4n)
 \le(A_\alpha+1)2^n.
\end{align*}

Finally, using $P_0=2^{\beta_d}$ and
\eqref{eq:step_size_iteration}, we obtain
\begin{align*}
 \log(16P_0^2)
 &=4\log2+2(\log2)\beta_d\\
 &\le
 4\log2+2(\log2)T_\alpha\bigl((4n)^d-1\bigr)\\
 &\le2(\log2)T_\alpha(4n)^d\\
 &\le\exp\bigl(T_\alpha+(A_\alpha+1)2^n\bigr)\\
 &\le\exp(B_\alpha2^n).
\end{align*}
Here we used $T_\alpha\ge2$ and
$\log(2(\log2)T_\alpha)\le T_\alpha$.
This proves~\eqref{eq:step_P_0_tau_0_size}.
\end{proof}

For the final choice of parameters, let $\Lambda_\alpha \geq 1$ be defined by
\begin{equation}\label{eq:step_K_alpha}
\Lambda_\alpha :=  \bigg\lceil(A_\alpha+1)\log2+\frac18\log(4B_\alpha)\bigg\rceil.
\end{equation}
Since $N\ge8$ and $n\le N(A_\alpha+1)+2$, we have
\[
 B_\alpha2^n
 \le4B_\alpha\exp\bigl((A_\alpha+1)(\log2)N\bigr)
 \le\exp(\Lambda_\alpha N).
\]
It follows from~\eqref{eq:step_P_0_tau_0_size} that
\begin{equation}\label{eq:step_height_bound}
 16P_0^2\le\exp\bigl(\exp(\exp(\Lambda_\alpha N))\bigr).
\end{equation}
In the symmetric case $\alpha=-1$, we may take $L_{-1}=10$.
Then
\[
 A_{-1}=\frac{1048576}{127},\qquad
 T_{-1}=20,\qquad B_{-1}=A_{-1}+21,
\]
and \[\Lambda_{-1}=5725.\]

\subsection{Subdivision estimates}
\label{subsec:step_subdivision}

We divide $[0,P_0)$ into small intervals whose subintervals will be used to construct the first characteristic function. For each positive integer $L$ and each $a=0,\ldots,LP_0-1$, set
\begin{equation}\label{eq:step_I_a_w_a}
 I_{a,L}:=[a/L,(a+1)/L),
 \qquad
 \tau_{a,L}:=\tau_0(a/L),
 \qquad
 w_{a,L}:=L\tau_{a,L}.
\end{equation}
For each fixed $L$, the intervals $I_{a,L}$ form a partition of $[0,P_0)$. Since $L$ is an integer, for every $k=0,\ldots,P_0-1$ we have
\[
 [k,k+1)
 =\bigcup_{a=kL}^{(k+1)L-1}
   I_{a,L}.
\]
Thus each $I_{a,L}$ is contained in an interval $[k,k+1)$ on which $\tau_0$ is constant, so $\tau_0$ has the constant value $\tau_{a,L}$ on $I_{a,L}$.

Applying Lemma~\ref{lem:step_absolute_convergence} to $g_0$, we obtain, for every $x\in X_0$ and every $f\in L^\infty(\R)$,
\[
 H_\alpha(f,g_0)(x)=\int_\R f(z)\,K_x(z)\,\rd z,
\]
where
\begin{equation}\label{eq:step_kernel_K_x(z)}
 K_x(z):=\frac{\operatorname{sgn}(\alpha)\,g_0(s_x(z))}{x-z},
 \qquad z\ne x,
\end{equation}
and $K_x(x):=0$, where $s_x(z)$ is as in~\eqref{eq:recur_s_x(z)}. In particular, for every bounded interval~$I$,
\begin{equation}\label{eq:step_interval_kernel}
 H_\alpha(\chi_I,g_0)(x)=\int_I K_x(z)\,\rd z,
 \qquad x\in X_0.
\end{equation}

Recall the set $D_{g,\tau}(x)$ defined
in~\eqref{eq:recur_D_g,tau}.
For fixed $x\in X_0$ and a positive integer $L$, call an interval $I_{a,L}$ \emph{exceptional} if its interior contains one of the points
\[
 (1-\alpha)x+\alpha j,\qquad j=0,\ldots,P_0,
\]
and define
\begin{equation}\label{eq:step_G_x,L}
\begin{split}
 G_{x,L}:=\Big\{a\in\{0,\ldots,LP_0-1\}&:\,
 I_{a,L}\text{ is nonexceptional}\\
 & \qquad \text{and }\big|I_{a,L}\cap D_{g_0,\tau_0}(x)\big|>0\Big\}.
\end{split}
\end{equation}

We note that the only possible discontinuities of
$g_0(s_x(z))$ occur at the points $(1-\alpha)x+\alpha j$, $j=0,\ldots,P_0$, and so this function is constant a.e.\ on each nonexceptional interval. For fixed $L$, the set $G_{x,L}$ consists of the indices of the nonexceptional intervals whose intersection with
$D_{g_0,\tau_0}(x)$ has positive measure.

The following lemma extends this information to a lower bound for $\tau_{a,L}K_x$ a.e.\ on each $I_{a,L}$ with $a\in G_{x,L}$ and gives the integral estimate needed for the passage to characteristic functions.

\begin{lemma}\label{lem:step_subdivision}
For every positive integer $L$, every $x\in X_0$,
and every $a\in G_{x,L}$, we have
\begin{equation}\label{eq:step_tau_a,L_K_x_bound}
 \tau_{a,L} \, K_x(z)>\frac12
 \qquad\text{for a.e.\ }z\in I_{a,L}.
\end{equation}
Moreover,
\begin{equation}\label{eq:step_subdivision_error}
 \min\!\bigg(1,\sum_{a\in G_{x,L}}\frac1{w_{a,L}}\bigg)
 \ge C_{g_0,\tau_0}(x)-\frac{P_0}{L}.
\end{equation}
Consequently, whenever $L\ge8P_0$,
\begin{equation}\label{eq:step_subdivision_bound}
 \int_{X_0}
 \min\!\bigg(1,\sum_{a\in G_{x,L}}\frac1{w_{a,L}}\bigg)\,\dx
 >\frac{P_0}{4}.
\end{equation}
\end{lemma}

\begin{proof}
Fix $x\in X_0$ and a positive integer $L$. We first prove~\eqref{eq:step_tau_a,L_K_x_bound}. Let $a\in G_{x,L}$. Since $|I_{a,L}\cap D_{g_0,\tau_0}(x)|>0$, choose
\[
 z_0\in \big(a/L,(a+1)/L\big)
       \cap D_{g_0,\tau_0}(x),
\]
and set $c:=g_0(s_x(z_0)) \geq 1$.
The interval $I_{a,L}$ is nonexceptional, so
$g_0(s_x(z))=c$ throughout its interior.
By the definition of $D_{g_0,\tau_0}(x)$,
\[s_x(z_0)<x,\qquad
 |z_0-x|<c \,\tau_{a,L}.
\]

We claim that $x$ does not belong to the interior of
$I_{a,L}$. Otherwise, since $s_x(x)=x$ and
$x\in X_0$, we would have
\[
 c=g_0(s_x(x))=g_0(x)=0,
\]
contradicting $c\ge1$. Thus $x-z_0$ has the same sign as $x-z$ for every $z$ in the interior of $I_{a,L}$.
Since
\[
 \frac{x-z_0}{\alpha}=x-s_x(z_0)>0,
\]
the definition of $K_x$ gives
\[
 K_x(z)=\frac{c}{|x-z|}
 \qquad\text{for a.e.\ }z\in I_{a,L}.
\]
For such $z$, using $L^{-1}\le1$ and $c \, \tau_{a,L}\ge2$, we obtain
\[
 |z-x|
 \le |z_0-x|+L^{-1}
 <c\tau_{a,L}+1
 \le2c\tau_{a,L}.
\]
Consequently,
\[
 \tau_{a,L}K_x(z)
 =\frac{c\tau_{a,L}}{|x-z|}
 >\frac12
 \qquad\text{for a.e.\ }z\in I_{a,L},
\]
which proves~\eqref{eq:step_tau_a,L_K_x_bound}.

We now prove~\eqref{eq:step_subdivision_error} by summing over all indices $a=0,\ldots,LP_0-1$.
For each such index, we have
\[
 \int_{I_{a,L}}\frac{\rd z}{\tau_0(z)}
 =\frac1{L\tau_{a,L}}
 =\frac1{w_{a,L}}.
\]
If $I_{a,L}$ is nonexceptional and $a\notin G_{x,L}$,
then $|I_{a,L}\cap D_{g_0,\tau_0}(x)|=0$, and hence
\[
 \int_{I_{a,L}\cap D_{g_0,\tau_0}(x)}
       \frac{\rd z}{\tau_0(z)}=0.
\]
Since the intervals $I_{a,L}$ partition $[0,P_0)$,
we obtain
\begin{align*}
 \int_{D_{g_0,\tau_0}(x)}\frac{\rd z}{\tau_0(z)}
 &=
 \sum_{a=0}^{LP_0-1}
 \int_{I_{a,L}\cap D_{g_0,\tau_0}(x)}
       \frac{\rd z}{\tau_0(z)}\\
 &\le
 \sum_{a\in G_{x,L}}\frac1{w_{a,L}}
 +
 \sum_{\substack{0\le a<LP_0\\
                 I_{a,L}\text{ exceptional}}}
 \int_{I_{a,L}}\frac{\rd z}{\tau_0(z)}\\
 &\le
 \sum_{a\in G_{x,L}}\frac1{w_{a,L}}
 +\frac{P_0+1}{2L}\\
 &\le
 \sum_{a\in G_{x,L}}\frac1{w_{a,L}}
 +\frac{P_0}{L}.
\end{align*}
Here we used $\tau_0\ge2$ and the fact that there are
at most $P_0+1$ exceptional intervals.
Taking the minimum with one, it follows that
\[
 \min\!\bigg(1,\sum_{a\in G_{x,L}}\frac1{w_{a,L}}\bigg)
 \ge C_{g_0,\tau_0}(x)-\frac{P_0}{L},
\]
which proves~\eqref{eq:step_subdivision_error}.

Finally, integrating over $X_0$, using $|X_0|\le P_0$ and Lemma~\ref{lem:step_Cgtau_bound}, we obtain
\begin{align*}
 \int_{X_0}
 \min\!\bigg(1,\sum_{a\in G_{x,L}}\frac1{w_{a,L}}\bigg)\,\dx
 &\ge
 \int_{X_0}C_{g_0,\tau_0}(x)\,\dx-\frac{P_0^2}{L}\\
 &>\frac{3P_0}{8}-\frac{P_0^2}{L}.
\end{align*}
If $L\ge8P_0$, then $P_0^2/L\le P_0/8$, which gives
\eqref{eq:step_subdivision_bound}.
\end{proof}

\section{Passage to characteristic functions}
\label{stripe:section}

In this section, we first approximate $g_0/U$ by a characteristic function, where $g_0$ is constructed in Section~\ref{sec:step} and $U$ is an integer satisfying the bound in Proposition~\ref{prop:char_extraction}
below. We then choose the first set as a union of short intervals, using the counting estimate from Lemma~\ref{lem:step_counting}.

\subsection{Lemmas for the second input}

The following estimate will control the change in the bilinear Hilbert transform when $g_0/U$ is replaced by a characteristic function.

\begin{lemma}\label{lem:char_integral_bound}
Let $\delta>0$, $D\ge0$, and $\varphi\in L^1(\R)$.
Suppose that $\varphi(y)=0$ for a.e.\ $y\in \R$ with $|y|<\delta$
and that
\[
 \left|\int_I \varphi(y)\,\rd y\right|\le D
\]
for every bounded interval $I$.
Then the integral below converges absolutely and
\[
 \left|\int_\R\frac{\varphi(y)}y\,\rd y\right|
 \le\frac{2D}{\delta}.
\]
\end{lemma}

\begin{proof}
We first note that the integral converges absolutely because
\[
 \int_\R\frac{|\varphi(y)|}{|y|}\,\rd y
 \le\frac{\|\varphi\|_1}{\delta}<\infty.
\]
Now, for $t\ge\delta$, set
\[
 F^\pm(t):=\int_\delta^t \varphi(\pm y)\,\rd y.
\]
Then $F^\pm(\delta)=0$, and the hypothesis gives
$|F^\pm(t)|\le D$, since $F^+(t)$ and $F^-(t)$ are the integrals of
$\varphi$ over $[\delta,t]$ and $[-t,-\delta]$, respectively.

For every $R>\delta$, the functions $F^\pm$ are absolutely
continuous on $[\delta,R]$, so integration by parts,
together with $F^\pm(\delta)=0$ and $|F^\pm|\le D$, gives
\begin{align*}
 \left|\int_\delta^R\frac{\varphi(\pm y)}y\,\rd y\right|
 &=\left|\frac{F^\pm(R)}R
       +\int_\delta^R\frac{F^\pm(y)}{y^2}\,\rd y\right|\\
 &\le\frac DR
       +D\left(\frac1\delta-\frac1R\right)
 =\frac D\delta.
\end{align*}
Since $\varphi$ vanishes a.e.\ on $(-\delta,\delta)$
and the integral converges absolutely, changing variables
on the negative half-line and applying these bounds yields
\begin{align*}
 \left|\int_\R\frac{\varphi(y)}y\,\rd y\right|=
 \lim_{R\to\infty}
 \left|\int_\delta^R
       \frac{\varphi(y)-\varphi(-y)}y\,\rd y\right|\le\frac{2D}{\delta},
\end{align*}
and concludes the proof.
\end{proof}

We next give a general approximation construction, which
will be applied to $g_0$ in
Proposition~\ref{prop:char_extraction}.

For the following two lemmas, let $P\ge1$ be an integer and let $g$ be a nonzero, nonnegative, integer-valued function vanishing outside $[0,P)$ and constant on each interval $[k,k+1)$ with $k\in\mathbb Z$.
Choose an integer \[U\ge\|g\|_\infty.\]
For every integer $m\ge1$, divide each interval $[k,k+1)$, $k=0,\ldots,P-1$, into $m$ intervals of length $1/m$ and retain the initial subinterval of length $g(k)/(Um)$ from each, defining
\begin{equation} \label{eq:char_B_m}
 B_m:=\bigcup_{k=0}^{P-1}\ \bigcup_{j=0}^{m-1}
 \left[k+\frac jm,\,
       k+\frac jm+\frac{g(k)}{Um}\right).
\end{equation}
By construction, it is clear that for each integer $m\ge1$, the set $B_m$ is contained in $[0,P) \cap \{g > 0\}$ and is a finite union of half-open intervals with endpoints in $(Um)^{-1}\mathbb Z$. By disjointness, we also have 
\[
 |B_m|=\sum_{k=0}^{P-1}m\,\frac{g(k)}{Um}
       =\frac{\|g\|_1}{U}.
\]

\begin{lemma}\label{lem:char_discrepancies}
For every integer $m\ge1$ and every bounded interval $I$,
\begin{equation}\label{eq:char_discrepancies}
 \bigg|\int_I
       \bigg(\chi_{B_m}(z)-\frac{g(z)}U\bigg)\,\rd z\bigg|
 \le\frac1{4m}.
\end{equation}
\end{lemma}

\begin{proof}
Set
\[
 F(z):=\int_{-\infty}^z
       \bigg(\chi_{B_m}(y)-\frac{g(y)}U\bigg)\,\rd y, \qquad z \in \R.
\]
For $k=0,\ldots,P-1$ and $j=0,\ldots,m-1$, the part of $[k+j/m,k+(j+1)/m)$ retained in $B_m$ has length $g(k)/(Um)$, while $g$ is constant and equal to $g(k)$ on this interval. Hence
\[
 \int_{k+j/m}^{k+(j+1)/m}
       \bigg(\chi_{B_m}(y)-\frac{g(y)}U\bigg)\,\rd y
 =\frac{g(k)}{Um}-\frac{g(k)}U\frac1m
 =0.
\]
Since both functions vanish on $(-\infty,0)$, we have
$F(0)=0$, and the preceding identities give
\[
 F\bigg(k+\frac{j+1}{m}\bigg)
 =F\bigg(k+\frac jm\bigg),
 \qquad
 k=0,\ldots,P-1,\quad j=0,\ldots,m-1.
\]
Applying these equalities successively, starting from
$F(0)=0$, we obtain
\[
 F\bigg(k+\frac jm\bigg)=0,
 \qquad
 k=0,\ldots,P-1,\quad j=0,\ldots,m-1.
\]
Thus, for $0\le t\le1/m$,
\[
 F\bigg(k+\frac jm+t\bigg)
 =
 \begin{cases}
  \displaystyle\bigg(1-\frac{g(k)}U\bigg)t,
    \qquad & \displaystyle 0\le t\le\frac{g(k)}{Um},\\[3mm]
  \displaystyle\frac{g(k)}U\bigg(\frac1m-t\bigg),
    \qquad & \displaystyle\frac{g(k)}{Um}\le t\le\frac1m.
 \end{cases}
\]
Both expressions are nonnegative and bounded above by
\[
 \frac1m\frac{g(k)}U\bigg(1-\frac{g(k)}U\bigg)
 \le\frac1{4m},
\]
since $0\le g(k)/U\le1$.
Together with the fact that $F$ is zero outside $[0,P)$, this gives
\[
 0\le F(z)\le\frac1{4m}
 \qquad\text{for every }z\in\R.
\]
For any bounded interval $I$ with endpoints $i_1\le i_2$, it follows that
\[
 \bigg|\int_I
       \bigg(\chi_{B_m}(z)-\frac{g(z)}U\bigg)\,\rd z\bigg|
 =|F(i_2)-F(i_1)|
 \le\frac1{4m},
\]
which proves~\eqref{eq:char_discrepancies}.
\end{proof}

We now use Lemma~\ref{lem:char_integral_bound} to estimate
the change in the transform when $g/U$ is replaced
by $\chi_{B_m}$.

\begin{lemma}\label{lem:char_stripe}
Let $\alpha\in\R\setminus\{0\}$, and suppose that
$E\subseteq\R$ and $\delta>0$ satisfy
\[
 \operatorname{dist}(x,\operatorname{supp}g)\ge\delta
 \qquad\text{for every }x\in E.
\]
Then, for every integer $m\ge1$, every bounded interval $I$, and every $x\in E$, the integrals $H_\alpha(\chi_I,\chi_{B_m})(x)$ and $H_\alpha(\chi_I,g)(x)$ converge absolutely and
\begin{equation}\label{eq:char_approx}
 \left|H_\alpha(\chi_I,\chi_{B_m})(x)
       -\frac1U H_\alpha(\chi_I,g)(x)\right|
 \le\frac1{2m\delta}.
\end{equation}
\end{lemma}

\begin{proof}
Fix an integer $m\ge1$, a bounded interval $I$, and
$x\in E$.
The construction in~\eqref{eq:char_B_m} gives
$B_m\subseteq\{g>0\}$ and hence
$\operatorname{supp}\chi_{B_m}\subseteq\operatorname{supp}g$,
so
\begin{equation}\label{eq:char_hole}
 \operatorname{dist}(x,\operatorname{supp}\chi_{B_m})
 \ge\delta.
\end{equation}
Set
\[
 \varphi_x(y):=\chi_I(x-\alpha y)
       \bigg(\chi_{B_m}(x-y)-\frac{g(x-y)}U\bigg).
\]
Since $g,\chi_{B_m}\in L^1(\R)$, we have
$\varphi_x\in L^1(\R)$, while the distance assumption
and~\eqref{eq:char_hole} imply that both $g(x-y)$ and
$\chi_{B_m}(x-y)$ vanish for a.e.\ $y \in \R$ with $|y|<\delta$. In particular, $\varphi_x(y)=0$ whenever $|y| < \delta$, and 
\begin{align*}
 \int_\R
 \frac{\chi_I(x-\alpha y)\,\chi_{B_m}(x-y)}{|y|}\,\rd y
 &\le\frac{|B_m|}{\delta}<\infty,\\
 \int_\R
 \frac{\chi_I(x-\alpha y)\,g(x-y)}{|y|}\,\rd y
 &\le\frac{\|g\|_1}{\delta}<\infty,
\end{align*}
so both transforms in~\eqref{eq:char_approx}
converge absolutely.

For every bounded interval $J$, the intersection
$J\cap\big\{y\in\R:\,x-\alpha y\in I\big\}$ is a bounded
interval or is empty.
Restricting the integral to this intersection, changing
variables $z=x-y$, and applying
\eqref{eq:char_discrepancies} therefore gives
\[
 \bigg|\int_J\varphi_x(y)\,\rd y\bigg|\le\frac1{4m}.
\]
Applying Lemma~\ref{lem:char_integral_bound} to $\varphi_x$
with $D=1/(4m)$ now yields
\begin{align*}
 \bigg|H_\alpha(\chi_I,\chi_{B_m})(x)
       -\frac1U H_\alpha(\chi_I,g)(x)\bigg|
 =\bigg|\int_\R\frac{\varphi_x(y)}y\,\rd y\bigg|\le\frac1{2m\delta},
\end{align*}
which proves~\eqref{eq:char_approx}.
\end{proof}

\subsection{Construction of the two sets}

We now use the functions constructed in Section~\ref{sec:step} to choose both sets.

\begin{proposition}\label{prop:char_extraction}
Fix an integer $N\ge8$, and let $P_0,g_0,\tau_0,X_0$ be as in Section~\ref{subsec:step_initial_functions}.
For every integer $U\ge8P_0^2$ and every $h_0>0$, there exist a number $h\in(0,h_0)$, sets $A,B\subseteq[0,P_0)$ of positive measure, each a finite union of half-open intervals of length~$h$, and a measurable set $E\subseteq X_0$ of positive measure such that, for every $x\in E$, the integral defining $H_\alpha(\chi_A,\chi_B)(x)$ converges absolutely and
\begin{equation}\label{eq:char_output}
 \big|H_\alpha(\chi_A,\chi_B)(x)\big|>\frac1{8U^2}.
\end{equation}
Moreover,
\begin{equation}\label{eq:char_quotient}
 \frac{|E|^2}{|A||B|}>2^{-8}U^2N.
\end{equation}
\end{proposition}

\begin{proof}
We divide the proof into four steps.

\smallskip
\emph{Step $1$. Construction of $B$.}
Set $L:=8P_0$ and recall the definitions of $I_{a,L}$,
$\tau_{a,L}$, and $w_{a,L}$ in~\eqref{eq:step_I_a_w_a}.
For $a=0,\ldots,LP_0-1$, write
\begin{equation}\label{eq:char_intervals_weights}
 I_a:=I_{a,L},\qquad
 \tau_a:=\tau_{a,L},\qquad
 w_a:=w_{a,L}=L\tau_a,
\end{equation}
and set
\begin{equation}\label{eq:char_short_intervals}
 A_a:=\big[a/L,a/L+\tau_a/U\big).
\end{equation}
Since $\|\tau_0\|_\infty\le P_0$ and $U\ge8P_0^2$,
we have $\tau_a/U\le1/L$. Thus $A_a\subseteq I_a$,
and the intervals $A_a$ are pairwise disjoint.

Choose an integer $m$ multiple of $L$ and sufficiently large that
\begin{equation}\label{eq:char_mesh_choice}
 m>4U^2,
 \qquad \frac1{Um}<h_0.
\end{equation}
Since $g_0$ is nonnegative and constant on each interval $[k,k+1)$ with $k\in\mathbb Z$, using~\eqref{eq:step_g1_bound}
and $U\ge8P_0^2$, we obtain
\[
 \|g_0\|_\infty\le\|g_0\|_1\le P_0/N\le U.
\]
Let $B:=B_m$, where $B_m$ is obtained by applying the construction in~\eqref{eq:char_B_m} to $g_0$ on~$[0,P_0)$.
Then
\begin{equation}\label{eq:char_B_mass}
 |B|=\frac{\|g_0\|_1}{U},
 \qquad
 \operatorname{supp}\chi_B\subseteq\operatorname{supp}g_0.
\end{equation}

For each $a=0,\ldots,LP_0-1$, define
\[
 \Phi_a(x):=U^2H_\alpha(\chi_{A_a},\chi_B)(x),
 \qquad x\in X_0.
\]
By Lemma~\ref{lem:char_stripe}, applied to the bounded interval $A_a$ and the set $X_0$ with $\delta=1/2$, the integral defining $H_\alpha(\chi_{A_a},\chi_B)(x)$
converges absolutely for every $x\in X_0$.
The functions $\Phi_a$ are measurable and real-valued, and~\eqref{eq:char_approx} together
with~\eqref{eq:char_mesh_choice} gives
\begin{equation}\label{eq:char_coefficient_error}
 \bigg|\Phi_a(x)-U H_\alpha(\chi_{A_a},g_0)(x)\bigg|
 \le\frac{U^2}{m}<\frac14,
 \qquad x\in X_0.
\end{equation}

\smallskip
\emph{Step $2$. Application of the subdivision lemma.}
Fix $x\in X_0$ and recall $G_{x,L}$ from~\eqref{eq:step_G_x,L}. For every $a\in G_{x,L}$, Lemma~\ref{lem:step_subdivision} gives
$\tau_a K_x(z)>1/2$ for a.e.\ $z\in I_a$.
Since $A_a\subseteq I_a$ and $|A_a|=\tau_a/U>0$,
the representation in~\eqref{eq:step_interval_kernel} gives
\[
 U H_\alpha(\chi_{A_a},g_0)(x)
 =U\int_{A_a}K_x(z)\,\rd z
 >\frac{U|A_a|}{2\tau_a}
 =\frac12.
\]
Together with~\eqref{eq:char_coefficient_error}, this implies $\Phi_a(x)>1/4$ for every $a\in G_{x,L}$. For each $x\in X_0$, define
\[
 J_x:=\big\{a\in\{0,\ldots,LP_0-1\}:\,
              \Phi_a(x)>1/4\big\}.
\]
The preceding argument gives $G_{x,L}\subseteq J_x$
for every $x\in X_0$. Using $w_a=L\tau_a$ and $L=8P_0$, Lemma~\ref{lem:step_subdivision} therefore yields
\begin{equation}\label{eq:char_capacity}
 \int_{X_0}\min\!\bigg(1,
       \sum_{a\in J_x}\frac1{w_a}\bigg)\,\dx
 \ge
 \int_{X_0}\min\!\bigg(1,
       \sum_{a\in G_{x,L}}\frac1{L\tau_a}\bigg)\,\dx
 >\frac{P_0}{4}.
\end{equation}

\smallskip
\emph{Step $3$. Counting estimates.}
Set
\[
 W:=\prod_{a=0}^{LP_0-1}\{0,\ldots,w_a-1\},
 \qquad K:=\#W,
\]
and for each $s=(s_0,\ldots,s_{LP_0-1})\in W$, define
\[
 A^{(s)}:=\bigcup_{\substack{0\le a<LP_0\\s_a=0}}A_a.
\]
Since the intervals $A_a$ are disjoint, absolute
convergence and linearity give
\[
 U^2H_\alpha(\chi_{A^{(s)}},\chi_B)(x)
 =\sum_{a=0}^{LP_0-1}\Phi_a(x)\chi_{\{0\}}(s_a),
 \qquad x\in X_0.
\]
For each $s\in W$, set
\begin{align*}
 E_s:=\big\{x\in X_0&:\,
       |H_\alpha(\chi_{A^{(s)}},\chi_B)(x)|>1/(8U^2)\big\},
\end{align*}
and
\[e_s:=|E_s|, \qquad 
  m_s:=|A^{(s)}|. \]
The sets $E_s$ are measurable, and for each fixed
$x\in X_0$, Lemma~\ref{lem:step_counting} applies to
the weights $w_a$ and coefficients $\Phi_a(x)$, with
$t=1/8$ and index set $J_x$, since
$\Phi_a(x)>1/4=2t$ for every $a\in J_x$.
Therefore,
\[
 \sum_{s\in W}\chi_{E_s}(x) = \# \big\{ s \in W : \, x \in E_s \big\}
 \ge\frac K4
 \min\!\bigg(1,\sum_{a\in J_x}\frac1{w_a}\bigg).
\]
Integrating over $X_0$ and using~\eqref{eq:char_capacity},
we obtain
\begin{equation}\label{eq:char_sum_e_s}
 \sum_{s\in W}e_s>\frac{KP_0}{16}.
\end{equation}

Moreover, the disjointness of the intervals $A_a$ gives
\[
 m_s=\sum_{a=0}^{LP_0-1}\frac{\tau_a}{U}
                    \chi_{\{0\}}(s_a).
\]
For each $a$, exactly $K/w_a$ elements $s\in W$ have
$s_a=0$. Summing over $s$, we therefore obtain
\begin{equation}\label{eq:char_sum_m_s}
 \sum_{s\in W}m_s
 =\sum_{a=0}^{LP_0-1}\frac{\tau_a}{U}\frac K{w_a}
 =\sum_{a=0}^{LP_0-1}\frac K{LU}
 =\frac{KP_0}{U}.
\end{equation}

\smallskip
\emph{Step $4$. Choice of $A$ and $E$.}
If $m_s=0$, then $A^{(s)}$ is empty and $e_s=0$.
Since $W$ is finite and the sum in~\eqref{eq:char_sum_m_s} is positive, the maximum
\[
 M:=\max_{\substack{s\in W\\m_s>0}}\frac{e_s^2}{m_s}
\]
exists, and the Cauchy--Schwarz inequality gives
\begin{align*}
 \bigg(\sum_{s\in W}e_s\bigg)^2
 &\le
 \bigg(\sum_{\substack{s\in W\\m_s>0}}\frac{e_s^2}{m_s}\bigg)
 \bigg(\sum_{s\in W}m_s\bigg)\\
 &\le KM\sum_{s\in W}m_s.
\end{align*}
Using~\eqref{eq:char_sum_e_s} and~\eqref{eq:char_sum_m_s}, we conclude that
\begin{equation} \label{eq:char_aux_M_bound}
 M\ge
 \frac{\Big(\sum_{s\in W}e_s\Big)^2}
      {K\sum_{s\in W}m_s}
 >\frac{(KP_0/16)^2}{K^2P_0/U}
 =\frac{UP_0}{2^8}.
\end{equation}

Choose $s\in W$ with $m_s>0$ attaining this maximum, and set $A:=A^{(s)}$ and $E:=E_s$.
Then $|A|>0$ and $|E|>0$, since $e_s^2=M m_s>0$,
while $|B|>0$ follows from~\eqref{eq:char_B_mass}.
The definition of $E_s$ gives~\eqref{eq:char_output}, and absolute convergence holds on $E\subseteq X_0$.
Moreover, using~\eqref{eq:char_aux_M_bound},
\[
 \frac{|E|^2}{|A||B|}
 =\frac{MU}{\|g_0\|_1}
 >\frac{U^2P_0}{2^8\|g_0\|_1}
 \ge2^{-8}U^2N,
\]
where the last inequality follows from
\eqref{eq:step_g1_bound}. This proves~\eqref{eq:char_quotient}.

Finally, the construction in~\eqref{eq:char_B_m} shows that the endpoints of the intervals defining~$B$ belong to $(Um)^{-1}\mathbb Z$.
Since $m$ is a multiple of $L$ and each $\tau_a$ is an integer, both endpoints of every $A_a$ also belong to $(Um)^{-1}\mathbb Z$.
Subdividing the intervals defining $A$ and $B$ therefore writes both sets as finite unions of half-open intervals of the common length $h:=1/(Um)<h_0$. This completes the proof.
\end{proof}

\section{Proof of Theorem~\ref{thm:main}}
\label{sec:proof}

In this section, we complete the proof of the main result by choosing the parameters in Proposition~\ref{prop:char_extraction}
and rescaling the resulting sets. The numerical estimates below are chosen for convenience and are not intended to be optimal.

\begin{proof}
By~\eqref{eq:inversion}, it suffices to prove the theorem
for $\alpha\in[-1,1)\setminus\{0\}$.
Let $\Lambda_\alpha\ge1$ be as in~\eqref{eq:step_K_alpha},
and set
\begin{equation}\label{eq:proof_lambda_alpha}
 \sigma_\alpha:=32\Lambda_\alpha,\qquad
 \lambda_\alpha:=\exp\big(\!-\exp(\exp\sigma_\alpha)\big).
\end{equation}
Fix $h_0>0$ and $0<\lambda<\lambda_\alpha$, and set
\begin{equation}\label{eq:proof_sigma_N}
 \sigma:=\log\log\log(1/\lambda),\qquad
 N:=\bigg\lfloor\frac{\sigma}{2\Lambda_\alpha}\bigg\rfloor.
\end{equation}
Since $\sigma>\sigma_\alpha$, we have
$\sigma/(2\Lambda_\alpha)>16$.
Using $\lfloor a\rfloor\ge a/2$ for $a\ge2$, we obtain
\begin{equation}\label{eq:proof_N_bound}
 N\ge\frac{\sigma}{4\Lambda_\alpha}\ge8.
\end{equation}
Construct $P_0,g_0,\tau_0,X_0$ as in
Section~\ref{subsec:step_initial_functions}
with this choice of $N$.
By~\eqref{eq:step_height_bound},
\[
 16P_0^2
 \le\exp\big(\exp(\exp(\Lambda_\alpha N))\big).
\]
Since $\Lambda_\alpha N\le\sigma/2$ and $e^\sigma-e^{\sigma/2}\ge\log4$ for $\sigma\ge2$,
it follows that
\[
 \log(16P_0^2)
 \le\exp(\exp(\sigma/2))
 \le\frac14\exp(\exp(\sigma))
 =\frac14\log(1/\lambda),
\]
and hence
\begin{equation}\label{eq:proof_P_0_bound}
 16P_0^2\le\lambda^{-1/4}.
\end{equation}

Choose
\[
 U:=\bigg\lfloor\sqrt{\frac1{8\lambda}}\bigg\rfloor.
\]
Since $\sigma_\alpha\ge32$, the choice of
$\lambda_\alpha$ and the inequality $\lambda<\lambda_\alpha$ give $\lambda\le2^{-10}$, so $\sqrt{1/(8\lambda)}\ge2$.
Using $a/2\le\lfloor a\rfloor\le a$ for $a\ge2$,
we obtain
\[
 \frac1{\sqrt{32\lambda}}
 \le U\le\frac1{\sqrt{8\lambda}}.
\]
Consequently, by~\eqref{eq:proof_P_0_bound},
\begin{equation}\label{eq:proof_U_size}
 U\ge\frac1{\sqrt{32\lambda}}
 \ge\lambda^{-1/4}
 \ge16P_0^2
 \ge8P_0^2,
\end{equation}
where the second inequality uses $\lambda\le2^{-10}$. We also have
\begin{equation}\label{eq:proof_lambda_U_bounds}
 \frac1{32}\le\lambda U^2\le\frac18.
\end{equation}

By~\eqref{eq:proof_U_size}, we may apply
Proposition~\ref{prop:char_extraction} with parameters
$U$ and $P_0h_0$ to infer the existence of a number
$\widetilde h\in(0,P_0h_0)$, sets
$\widetilde A,\widetilde B\subseteq[0,P_0)$ of positive
measure, each a finite union of half-open intervals of
length $\widetilde h$, and a measurable set $E\subseteq X_0$
of positive measure such that
\[
 \big|H_\alpha(\chi_{\widetilde A},\chi_{\widetilde B})(x)\big|
 >\frac1{8U^2}\ge\lambda
 \qquad\text{for every }x\in E,
\]
and
\[
 \frac{|E|^2}{|\widetilde A||\widetilde B|}
 >2^{-8}U^2N,
\]
where the inequality $1/(8U^2)\ge\lambda$ follows from~\eqref{eq:proof_lambda_U_bounds}.

The integral defining
$H_\alpha(\chi_{\widetilde A},\chi_{\widetilde B})(x)$
converges absolutely for every $x\in E$.
Since $\lambda U^2\ge2^{-5}$,
we obtain
\begin{equation}\label{eq:characteristic-ratio}
 \frac{\lambda|E|^2}{|\widetilde A||\widetilde B|}
 >2^{-13}N.
\end{equation}

Set
\[
 A:=P_0^{-1}\widetilde A,\qquad
 B:=P_0^{-1}\widetilde B,\qquad
 E':=P_0^{-1}E,\qquad
 h:=\frac{\widetilde h}{P_0}.
\]
Then $A,B\subseteq[0,1)$ have positive measure and
are finite unions of half-open intervals of the common
length $h\in(0,h_0)$.
For every $x\in E$, the change of variables $v=P_0y$
gives
\[
 H_\alpha(\chi_A,\chi_B)(x/P_0)
 =H_\alpha(\chi_{\widetilde A},\chi_{\widetilde B})(x),
\]
with absolute convergence on both sides.
Thus $\big|H_\alpha(\chi_A,\chi_B)\big|>\lambda$ on $E'$, and so
\[E' \subseteq \big\{x\in\R:\,
 |H_\alpha(\chi_A,\chi_B)(x)|>\lambda\big\}.\]
Moreover,
\[
 |A|=\frac{|\widetilde A|}{P_0},\qquad
 |B|=\frac{|\widetilde B|}{P_0},\qquad
 |E'|=\frac{|E|}{P_0}.
\]
Using these identities, \eqref{eq:characteristic-ratio},
and~\eqref{eq:proof_N_bound}, we obtain
\begin{align*}
 \frac{\lambda
 \big|\big\{x\in\R:\,
 |H_\alpha(\chi_A,\chi_B)(x)|>\lambda\big\}\big|^2}
 {|A||B|}
 &\ge\frac{\lambda|E'|^2}{|A||B|}\\
 &>2^{-13}N
 \ge\frac{2^{-15}}{\Lambda_\alpha}\,\sigma.
\end{align*}
Since $\sigma=\log\log\log(1/\lambda)$, this proves the theorem with
\begin{equation}\label{eq:proof_c_alpha}
 c_\alpha:=\frac{2^{-15}}{\Lambda_\alpha}.
\end{equation}
In particular, using $\Lambda_{-1}=5725$
we obtain
\[
 c_{-1}
 =\frac1{2^{15}\cdot5725}
 =\frac1{187596800}
 \approx5.3305813\times10^{-9}.
\]
\end{proof}

\section*{Acknowledgments}
The authors acknowledge the use of AI tools to assist with some constructions in this work. The authors developed the mathematical ideas and take full responsibility for the content and correctness of the manuscript.

%The research of Alves was supported by King Abdullah University of Science and Technology (KAUST) under Award No.\ ORFS-CRG12-2024-6430. The research of Grafakos...

\end{document}